\documentclass[11pt]{article}
\usepackage{cite}
\usepackage{mathrsfs}
\usepackage{amssymb,amsfonts,amsmath}
\usepackage{enumerate}
\usepackage{indentfirst}
\usepackage{latexsym,bm}
\usepackage[noend]{algpseudocode}
\usepackage{algorithmicx,algorithm}
\usepackage{algorithm}
\usepackage{makeidx}
\usepackage{fancybox}
\usepackage{color}
\usepackage{multicol}
\usepackage[latin1]{inputenc}
\usepackage{graphicx}
\usepackage{epstopdf}
\usepackage{pstricks,pst-node,pst-text,pst-3d}
\usepackage{tikz}
\usepackage{makecell}
\usepackage{caption}
\usepackage[urlcolor=blue]{hyperref}
\newtheorem{thm}{Theorem}[section]
\newtheorem{cor}[thm]{Corollary}
\newtheorem{lem}[thm]{Lemma}
\newtheorem{op}{Problem}[section]

\newtheorem{pro}[thm]{Proposition}

\newtheorem{conj}{Conjecture}[section]

\newenvironment{pf}{{\noindent \it \bf Proof:}}{{\hfill$\Box$}\\}

\begin{document}

\title{\bf Steiner Type Packing Problems in Digraphs: A Survey}
\author{
Yuefang Sun\\
School of Mathematics and Statistics, Ningbo University,\\
Ningbo 315211, P. R. China\\ 
Email address: sunyuefang@nbu.edu.cn}
\maketitle

\tableofcontents
\newpage
\begin{abstract}

Graph packing problem is one of the central problems in graph theory and combinatorial optimization. The famous Steiner tree packing problem in undirected graphs has become an well-established area. It is natural to extend this problem to digraphs, and such problems in digraphs are called directed Steiner type packing problems.

In this survey we overview known results on several directed Steiner type packing problems. The paper is divided into seven sections: introduction, directed Steiner tree packing problem, directed Steiner path packing problem, directed  pendant Steiner tree packing problem, strong subgraph packing problem, strong arc decomposition problem, directed Steiner cycle packing problem. This survey also contains some conjectures and open problems for further study.
\vspace{0.3cm}
\vspace{0.3cm}
\end{abstract}
\newpage

\section{Introduction}\label{sec:intro}

In this section, we introduce the backgrounds of undirected Steiner tree packing problem, and several directed Steiner type packing problems. We also introduce terminology and notation which appear in the survey. 

\subsection{Backgrounds}

We refer the readers to \cite{Bang-Jensen-Gutin} for graph theoretical notation and terminology not given here. Throughout this survey, unless otherwise stated, paths and cycles are always assumed to be directed, and all digraphs considered in this paper have no parallel arcs or loops. We use $[n]$ to denote the set of all natural numbers from 1 to $n$.

\subsubsection{Steiner tree packing problem}

For a graph $G=(V,E)$ and a (terminal) set $S\subseteq V$ of at
least two vertices, an {\em $S$-Steiner tree} or, simply, an {\em
$S$-tree} is a tree $T$ of $G$ with $S\subseteq V(T)$. By definition, the concept of Steiner tree is a common generalization of the concepts of path and spanning tree. Two $S$-trees are said to be {\em 
edge-disjoint} if they have no common edge. Two edge-disjoint $S$-trees are said to be {\em internally-disjoint} if the set of common vertices of them is exactly $S$. The basic problem of {\sc Steiner Tree
Packing} is defined as follows: 

\vspace{2mm}

\noindent {\bf Steiner Tree Packing (STP):}
The input consists of an undirected graph $G$ and a subset of vertices
$S\subseteq V(G)$, the goal is to find a largest collection of pairwise edge-disjoint $S$-trees.

\vspace{2mm}

From a theoretical perspective, both extremes of this problem are fundamental
theorems in combinatorics. One extreme of the problem is when $|S|=2$, in this case, finding edge-disjoint Steiner trees is equivalent to finding edge-disjoint paths between the two terminals in $S$, and so the problem becomes
the well-known Menger's theorem. The other extreme is when $|S|=n$, in this case, edge-disjoint Steiner trees are just edge-disjoint spanning trees of $G$, and so the problem becomes the classical Nash-Williams--Tutte theorem \cite{Nash-Williams, Tutte}.

From a practical perspective, the Steiner tree packing problem has applications in VLSI circuit design \cite{Grotschel-Martin-Weismantel, Sherwani}. 
In this application, a Steiner tree is needed to share an electronic signal by a set of terminal nodes. Another application arises in the Internet Domain \cite{Li-Mao5}: Let a given graph $G$ represent a network, we choose arbitrary $k$ vertices as nodes such that one of them is a {\em broadcaster}, and all other nodes are either {\em users} or {\em routers} (also called {\em switches}). The broadcaster wants to broadcast as many streams of movies as possible, so that the users have the maximum number of choices. Each stream of movie is broadcasted via a tree connecting all the users and the broadcaster. Hence we need to find the maximum number Steiner trees connecting all the users and the broadcaster, and it is a Steiner tree packing problem.

The Steiner tree packing problem (in undirected graphs) has attracted much attention from researchers in the area of graph theory, combinatorial optimization and theoretical computer sciences, and has become a well-established research topic \cite{Cheriyan-Salavatipour, DeVos-McDonald-Pivotto, Frank-Kiraly-Kriesell, Grotschel-Martin-Weismantel2, Grotschel-Martin-Weismantel3, Grotschel-Martin-Weismantel4, Grotschel-Martin-Weismantel5, Grotschel-Martin-Weismantel, Jain-Mahdian-Salavatipour, Jeong-Lee-Park-Park, Kriesell, Kriesell2, Kriesell3, Lau, Li-Mao5, Pulleyblank, Sherwani, Uchoa, West-Wu}.

\subsubsection{Directed Steiner tree packing problem}

An {\em out-tree} (resp. {\em in-tree}) is an oriented tree in which every vertex
except one, called the {\em root}, has in-degree (resp.  out-degree) one.
An {\em out-branching} (resp. {\em in-branching}) of $D$ is a spanning out-tree
(resp.  in-tree) in $D$. For a digraph $D=(V(D), A(D))$, and a set $S\subseteq V(D)$ with
$r\in S$ and $|S|\geq 2$, a {\em directed $(S, r)$-Steiner tree} or,
simply, an {\em $(S, r)$-tree} is an out-tree $T$ rooted at $r$ with
$S\subseteq V(T)$ \cite{Cheriyan-Salavatipour}. Two $(S, r)$-trees are said to be {\em arc-disjoint} if they have no common arc. Two arc-disjoint $(S, r)$-trees are said to be {\em internally-disjoint} if the set of common vertices of them is exactly $S$. Let $\kappa_{S,r}(D)$ (resp. $\lambda_{S,r}(D)$) be the maximum
number of pairwise internally-disjoint (resp. arc-disjoint) $(S, r)$-trees in $D$.

Cheriyan and Salavatipour \cite{Cheriyan-Salavatipour} introduced
and studied the following two directed Steiner tree packing
problems:

\vspace{2mm}

\noindent {\bf Arc-disjoint Directed Steiner Tree Packing (ADSTP):}
The input consists of a digraph $D$ and a subset of vertices
$S\subseteq V(D)$ with a root $r$, the goal is to find a largest collection of pairwise arc-disjoint $(S, r)$-trees.

\vspace{2mm}

\noindent {\bf Internally-disjoint Directed Steiner Tree Packing (IDSTP):}
The input consists of a digraph $D$ and a subset of vertices $S\subseteq V(D)$ with a root $r$, the goal is to find a largest collection of pairwise internally-disjoint $(S, r)$-trees.

\vspace{2mm}

It is worth mentioning that directed Steiner tree packing problems are related to two important theorems in digraph theory: when $|S|=2$, that is, $S=\{r, x\}\subseteq V(D)$, then $\lambda_{S,r}(D)=\lambda_D(r,x)$ (resp. $\kappa_{S,r}(D)=\kappa_D(r,x)$), the local arc-strong (resp. vertex-strong) connectivity of $r, x$ in $D$, and so now they are related to Menger's theorem; when $|S|=n$, directed Steiner tree packing problems are equivalent to finding a largest collection of pairwise arc-disjoint out-branchings rooted at $r$, and therefore are related to the famous Edmonds' Branching Theorem (Theorem~\ref{thm2-19}).

\subsubsection{Directed Steiner path packing problem}

Sun and Zhang~\cite{Sun5} introduced the concept of directed Steiner path packing which could be seen as a restriction of the directed Steiner tree packing problem. For a digraph $D=(V(D), A(D))$, and a set $S\subseteq V(D)$ with $r\in S$ and $|S|\geq 2$, a {\em directed $(S, r)$-Steiner path} or, simply, an {\em $(S, r)$-path} is a directed path $P$ started at $r$ with $S\subseteq V(P)$. Observe that the directed Steiner path is a directed Steiner tree which is a directed path, and it is also a generalization of the directed Hamiltonian path.

Two $(S, r)$-paths are said to be {\em arc-disjoint} if they have no common arc. Two arc-disjoint $(S, r)$-paths are said to be {\em internally-disjoint} if the set of common vertices of them is exactly $S$. Let $\kappa^p_{S,r}(D)$ (resp. $\lambda^p_{S,r}(D)$) be the maximum number of pairwise internally-disjoint (resp. arc-disjoint) $(S, r)$-paths in $D$. The directed Steiner path packing problems can be defined as follows:

\vspace{2mm}

\noindent {\bf Arc-disjoint Directed Steiner Path Packing (ADSPP):}
The input consists of a digraph $D$ and a subset of vertices
$S\subseteq V(D)$ with a root $r$, the goal is to find a largest
collection of pairwise arc-disjoint $(S, r)$-paths.

\vspace{2mm}

\noindent {\bf Internally-disjoint Directed Steiner Path Packing (IDSPP):}
The input consists of a
digraph $D$ and a subset of vertices $S\subseteq V(D)$ with a root
$r$, the goal is to find a largest collection of pairwise internally-disjoint
$(S, r)$-paths.

\subsubsection{Directed pendant Steiner Tree packing problem}

Yu and Sun~\cite{Yu-Sun2} introduced another restriction of the directed Steiner tree packing problem. If each vertex of $S$ has degree one in an $(S,r)$-tree $T$, then $T$ is called a {\em pendant $(S,r)$-tree}. 
Observe that in a pendant $(S, r)$-tree,  the in-degree of $r$ and the out-degree of each vertex $v\in S\setminus \{r\}$ is zero, and the out degree of $r$ and the in-degree of each vertex $v\in S\setminus \{r\}$ is one. Two pendant $(S, r)$-trees are said to be {\em arc-disjoint} if they have no common arc. Two arc-disjoint pendant $(S, r)$-trees are said to be {\em internally-disjoint} if the set of common vertices of them is exactly $S$. Let $\tau_{S,r}(D)$ (resp. $\tau'_{S,r}(D)$) be the maximum number of pairwise internally-disjoint (resp. arc-disjoint) pendant $(S, r)$-trees in $D$.

Yu and Sun~\cite{Yu-Sun2} introduced and studied the following problem:

\vspace{2mm}

\noindent {\bf Internally-disjoint Directed Pendant Steiner Tree Packing (IDPSTP):}
The input consists of a digraph $D$ and a subset of vertices $S\subseteq V(D)$ with a root $r$, the goal is to find a largest collection of pairwise internally-disjoint pendant $(S, r)$-trees.

\vspace{2mm}

Similarly, we can define the following problem:

\vspace{2mm}

\noindent {\bf Arc-disjoint Directed Pendant Steiner Tree Packing (ADPSTP):}
The input consists of a digraph $D$ and a subset of vertices
$S\subseteq V(D)$ with a root $r$, the goal is to find a largest collection of pairwise arc-disjoint pendant $(S, r)$-trees.

\subsubsection{Strong arc decomposition problem}

The following famous Edmonds' Branching Theorem is a fundamental theorem in the area of digraph packing theory.

\begin{thm}(Edmonds' Branching Theorem)\label{thm2-19}\cite{Edmonds}
A digraph $D = (V, A)$ with a special vertex $s$ has $k$ pairwise arc-disjoint out-branchings rooted at $s$ if and only if there are $k$ arc-disjoint $(s,v)$-paths in $D$ for every $v \in V-s$.
\end{thm}

Furthermore, there exists a polynomial algorithm for finding $k$ pairwise arc-disjoint out-branchings from a given root $s$ if they exist. However, if we ask for the existence of a pair of arc-disjoint branchings $B^+_s$, $B^-_s$ such that the former is an out-branching rooted at $s$ and the latter is an in-branching rooted at $s$, then the problem becomes NP-complete (see Section~9.6 of \cite{Bang-Jensen-Gutin}). In connection with this problem, Thomassen  posed the following conjecture:

\begin{conj}\label{conj4-1}\cite{Thomassen}
There exists an integer $N$ so that every $N$-arc-strong digraph $D$ contains a pair of arc-disjoint in- and out-branchings.
\end{conj}


A digraph $D=(V,A)$ has a {\em strong arc decomposition} if $A$ has two disjoint sets $A_1$ and $A_2$ such that both
$(V,A_1)$ and $(V,A_2)$ are strong \cite{Bang-Jensen-Gutin-Yeo, Bang-Jensen-Yeo}. Bang-Jensen and Yeo generalized Conjecture~\ref{conj4-1} as follows.\footnote{Every strong digraph $D$ has an out- and in-branching rooted at any vertex of $D$.}

\begin{conj}\label{conj4-2}\cite{Bang-Jensen-Yeo}
There exists an integer $N$ so that every $N$-arc-strong digraph $D$ has a strong arc decomposition.
\end{conj}

For a general digraph $D$, it is a hard problem to decide whether $D$ has a decomposition into two strong spanning subdigraphs.

\begin{thm}\label{thm4-13}\cite{Bang-Jensen-Yeo}
It is NP-complete to decide whether a digraph has a strong arc decomposition.
\end{thm}

\subsubsection{Strong subgraph packing problem}

There is another way to extend the Steiner tree packing problem to digraphs, note that an $S$-Steiner tree is a connected
subgraph of $G$ containing $S$. In fact, in the definition of Steiner tree packing problem, we could replace ``an $S$-Steiner tree'' by ``a connected subgraph of $G$ containing $S$". Therefore, we define the strong subgraph packing problem by replacing ``connected'' with ``strongly connected'' (or, simply, ``strong'') as follows. A digraph $D$ is {\em strongly connected} (or, {\em strong}), if for any pair of vertices $x, y\in V(D)$, there is a path from $x$ to $y$ in $D$, and vice versa. Let $D=(V(D),A(D))$ be a digraph of order $n$, $S\subseteq V$ a $k$-subset of $V(D)$ and $2\le k\leq n$. A strong subgraph $H$ of $D$ is called an {\em $S$-strong
subgraph} if $S\subseteq V(H)$. Two $S$-strong subgraphs are said to be {\em arc-disjoint} if they have no common arc. Furthermore, two arc-disjoint $S$-strong subgraphs are said to be {\em internally-disjoint} if the set of common vertices of them is exactly $S$. We use $\kappa_S(D)$~(resp. $\lambda_S(D)$) to denote the maximum number of pairwise internally-disjoint (resp. arc-disjoint) $S$-strong subgraphs in $D$ \cite{Sun-Gutin2, Sun-Gutin-Yeo-Zhang}.


Sun, Gutin and Zhang~\cite{Sun-Zhang} introduced the following two types of strong subgraph packing problems in digraphs which are analogs of directed Steiner tree packing problems:

\vspace{2mm}

\noindent {\bf Arc-disjoint strong subgraph packing (ASSP):}
The input consists of a digraph $D$ and a subset of vertices $S\subseteq V(D)$, the goal is to find a largest collection of pairwise  arc-disjoint $S$-strong subgraphs.

\vspace{2mm}

\noindent {\bf Internally-disjoint strong subgraph packing (ISSP):}
The input consists of a digraph $D$ and a subset of vertices $S\subseteq V(D)$, the goal is to find a largest collection of pairwise internally-disjoint $S$-strong subgraphs.

\vspace{2mm}

Observe that a digraph $D$ has a strong arc decomposition if and only if $\lambda_{V(D)}(D)\geq 2$ (or, $\kappa_{V(D)}(D)\geq 2$). Therefore, the problem of {\sc ASSP} (or {\sc ISSP}) could be seen as an extension of the strong arc decomposition problem.

\subsubsection{Directed Steiner cycle packing problem}

Let $D=(V(D),A(D))$ be a digraph of order $n$, $S\subseteq V$ a $k$-subset of $V(D)$ and $2\le k\leq n$. A directed cycle $C$ of $D$ is called a {\em directed $S$-Steiner cycle} or, simply, an {\em $S$-cycle} if $S\subseteq V(C)$. It is worth noting that Steiner cycles have applications in the optimal design of reliable telecommunication and transportation networks \cite{Steiglitz-Weiner-Kleitman}. 
Two $S$-cycles are said to be {\em arc-disjoint} if they have no common arc. Furthermore, two arc-disjoint $S$-cycles are {\em internally-disjoint} if the set of common vertices of them is exactly $S$. We use $\kappa^c_S(D)$~(resp. $\lambda^c_S(D)$) to denote the maximum number of pairwise internally-disjoint (resp. arc-disjoint) $S$-cycles in $D$ \cite{Sun4}.

Sun and Jin~\cite{Sun4} defined the following two directed Steiner cycle packing problems:

\vspace{2mm}

\noindent {\bf Arc-disjoint directed Steiner cycle packing (ADSCP):}
The input consists of a digraph $D$ and a subset of vertices $S\subseteq V(D)$, the goal is to find a largest collection of pairwise  arc-disjoint $S$-cycles.

\vspace{2mm}

\noindent {\bf Internally-disjoint directed Steiner cycle packing (IDSCP):}
The input consists of a digraph $D$ and a subset of vertices $S\subseteq V(D)$, the goal is to find a largest collection of pairwise internally-disjoint $S$-cycles.

\vspace{2mm}

By definition, directed Steiner cycle packing problems are restrictions of  strong subgraph packing problems, as an $S$-cycle is also an $S$-strong subgraph. However, there are differences between them. For example, observe that it can be decided in polynomial-time whether $\kappa_S(D)\geq 1$, but it is NP-complete to decide whether $\kappa^c_S(D)\geq 1$ even restricted to an Eulerian digraph $D$ \cite{Sun4}. 

Directed Steiner cycle packing problems 
are also related to other problems in graph theory. When $|S|=n$, an $S$-cycle is a Hamiltonian cycle. Therefore, directed Steiner cycle packing problems generalize the Hamiltonian cycle packing problem (and therefore the Hamiltonian decomposition problem). A digraph $D$ is {\em $k$-cyclic} if it has a cycle containing the vertices  $x_1, x_2, \dots, x_k$ for every choice of $k$ vertices. Note that the notion of $k$-cyclic attracts the attention of some researchers, such as \cite{Kuhn2008}. By definition, a digraph is $k$-cyclic if and only if $\kappa^c_S(D)\geq 1$ (resp. $\lambda^c_S(D)\geq 1$) for every $k$-subset $S$ of $V(D)$.


\subsubsection{Directed Steiner type packing problem}

Generally, for a digraph $D$ and a subset of vertices $S\subseteq V(D)$, a {\em directed $S$-Steiner type subgraph} is a subgraph $H$ of $D$ such that $S\subseteq V(H)$.
Two directed $S$-Steiner type subgraphs are said to be {\em arc-disjoint} if they have no common arc. Furthermore, two arc-disjoint directed $S$-Steiner type subgraphs are said to be {\em internally-disjoint} if the set of common vertices of them is exactly $S$. We now define the following problems:

\vspace{2mm}

\noindent {\bf Arc-disjoint directed Steiner type subgraph packing (ADSTSP):}
The input consists of a digraph $D$ and a subset of vertices $S\subseteq V(D)$, the goal is to find a largest collection of pairwise  arc-disjoint directed $S$-Steiner type subgraphs.

\vspace{2mm}

\noindent {\bf Internally-disjoint directed Steiner type subgraph packing (IDSTSP):}
The input consists of a digraph $D$ and a subset of vertices $S\subseteq V(D)$, the goal is to find a largest collection of pairwise internally-disjoint directed  $S$-Steiner type subgraphs.

\vspace{2mm}

\begin{figure}[!hbpt]
		\begin{center}
            \includegraphics[scale=0.3]{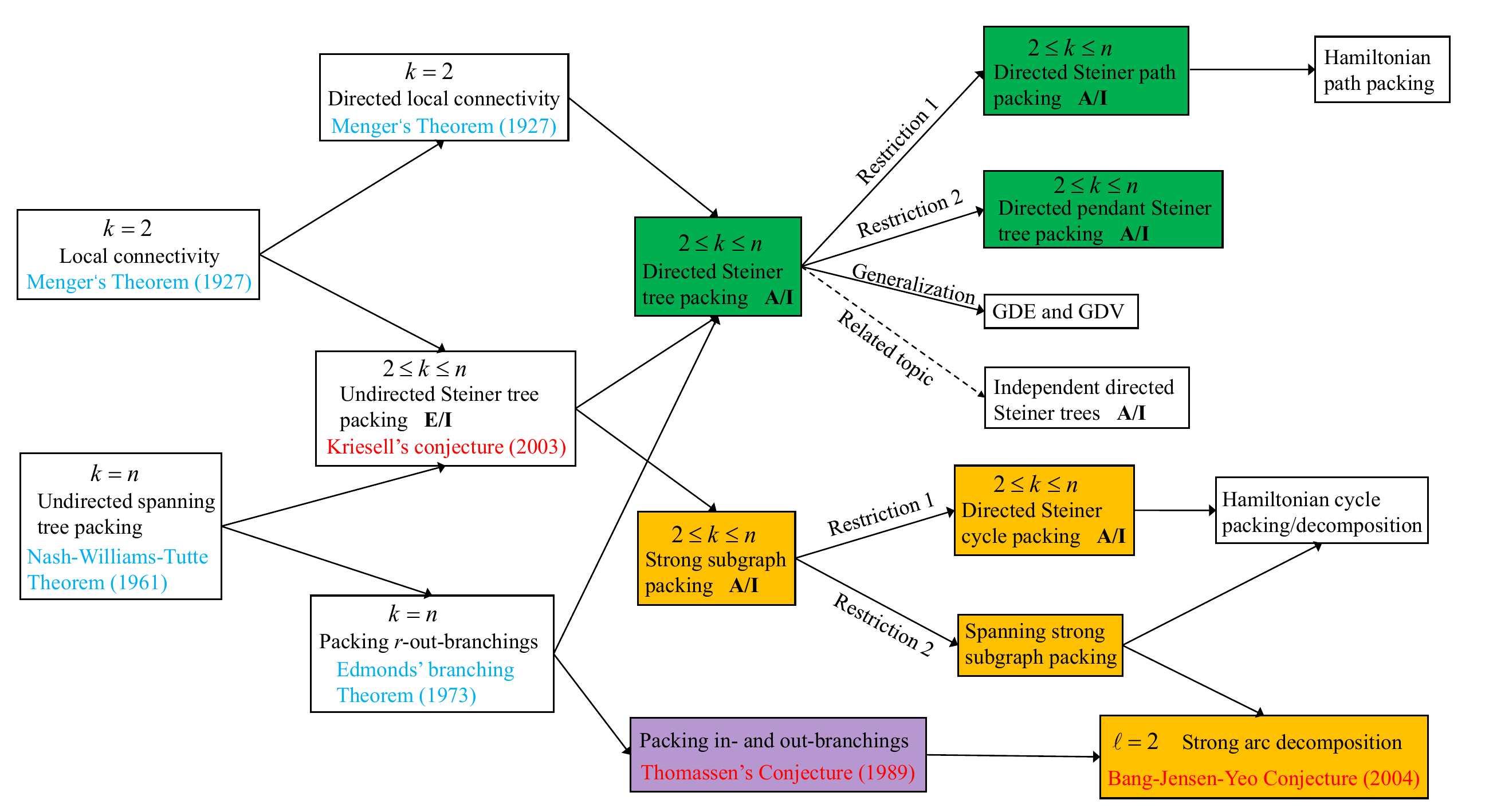}
		\end{center}
		\caption{Directed Steiner type packing and related topics.}\label{figure0}
	\end{figure}

{\sc ADSTSP} and {\sc IDSTSP} are also called {\em directed Steiner type packing problems}. By definition, all of {\sc ADSTP}, {\sc IDSTP}, {\sc ADSPP}, {\sc IDSPP}, {\sc ADPSTP}, {\sc IDPSTP}, {\sc ASSP}, {\sc ISSP} (and therefore strong arc decomposition problem), {\sc ADSCP} and {\sc IDSCP} belong to this type of problem. There are some papers on Steiner type packing problems in digraphs, such as \cite{Bai-Sun-Wang-Yu, Bang-Jensen-Gutin-Yeo, Bang-Jensen-Huang2012, Bang-Jensen-Wang2025, Bang-Jensen-Yeo, Cheriyan-Salavatipour, Dong-Gutin-Sun, Sun, Sun-Gutin2, Sun4, Sun5, Sun-Gutin2, Sun-Gutin-Ai, Sun-Gutin-Yeo-Zhang, Sun-Jin, Sun-Yeo, Sun-Zhang, Wang-Sun, Yu-Sun1, Yu-Sun2, Yu-Sun3}. The readers can see \cite{Sun-book} for a monograph on this topic. 
In addition, the relationship diagram between directed Steiner type packing and related topics is shown in Figure~\ref{figure0}.

Nowadays, an increasing number of researchers are conducting research in this topic, and new results keep emerging. 
In this survey, we try to collect most of known results on directed Steiner type packing problems. After an introductory section, the paper will be divided into six sections: directed Steiner tree packing problem, directed Steiner path packing problem, directed pendant Steiner tree packing problem, strong subgraph packing problem, strong arc decomposition problem, directed Steiner cycle packing problem. This survey also contains some conjectures and open problems for further study.

\subsection{Further terminology and notation}
A digraph is {\em connected} if its underlying graph is connected. A digraph $D$ is {\em semicomplete} if for every distinct $x,y\in V(D)$ at least one of the arcs $xy,yx$ is in $D$. A digraph $D$ is {\em locally semicomplete} if $N^-(x)$ and $N^+(x)$ induce semicomplete digraphs for every vertex $x$ of $D$. A digraph is {\em locally in-semicomplete} (resp. {\em locally out-semicomplete}) if $N^-(x)$ (resp. $N^+(x)$) induces a semicomplete digraph for every vertex $x$ of $D$. A {\em split digraph} is a digraph whose vertex set is a disjoint union of two nonempty sets $V_1$ and $V_2$ such that $V_1$ is an independent set and the subdigraph induced by $V_2$ is semicomplete. A digraph $D$ is {\em quasi-transitive}, if for any triple $x,y,z$ of distinct vertices of $D$, if $xy$ and $yz$ are arcs of $D$ then either $xz$ or $zx$ or both are arcs of $D.$  A digraph $D$ is {\em symmetric} if $xy \in A(D)$ then $yx \in A(D)$. In other words, a symmetric digraph $D$ can be obtained from its underlying undirected graph $G$ by replacing each edge of $G$ with the corresponding arcs of both directions, that is, $D=\overleftrightarrow{G}.$ For a digraph $D$, its {\em reverse} $D^{\rm rev}$ is a digraph with the same vertex set such that $xy\in A(D^{\rm rev})$ if and only if $yx\in A(D)$. Note that if a digraph $D$ is {\em symmetric} then $D^{\rm rev}=D$. A digraph is {\em $k$-regular} if the out-degree and in-degree of every vertex equals $k$. A digraph $D$ is {\em Eulerian} if $D$ is connected and $d^+(x)=d^-(x)$ for every vertex $x\in V(D)$. Observe that a connected symmetric digraph is Eulerian.

Let $T$ be a digraph with vertices $u_1, \dots, u_t$ ($t\ge 2$) and let $H_1, \dots, H_t$ be digraphs such that $H_i$ has vertices $u_{i,j_i}$, where $j_i\in [n_i].$  Let $n_0=\min\{n_i\mid i \in [t]\}$. The {\em composition} $Q=T[H_1, \dots, H_t]$ is a digraph with vertex set $$V(Q)=\{u_{i,j_i}\mid i \in [t], j_i\in [n_i]\}$$ and arc set $$A(Q)=\cup^t_{i=1}A(H_i)\cup \{u_{i,j_i}u_{p,q_p}\mid u_iu_p\in A(T), j_i\in [n_i], q_p\in [n_p]\}.$$ The composition $Q=T[H_1, \dots, H_t]$ is {\em proper} if $Q\neq T$, i.e. at least one $H_i$ is nontrivial, and is {\em semicomplete} (resp. {\em symmetric}) if $T$ is semicomplete (resp. symmetric). If $Q=T[H_1, \dots, H_t]$ and none of the digraphs $H_1, \dots, H_t$ has an arc, then $Q$ is an {\em extension} of $T$. For any set $\Phi$ of digraphs, $\Phi^{ext}$ denotes the (infinite) set of all extensions of digraphs in $\Phi$, which are called {\em extended $\Phi$-digraphs}. For example, if $\Phi$ is the set of all semicomplete digraphs, then $\Phi^{ext}$ denotes the set of all extended semicomplete digraphs.  The readers can see \cite{Sun2} for a recent survey on the topic of digraph compositions. 

We now introduce the definitions of several product digraphs \cite{Hammack}. The {\em Cartesian product} $G\square H$ of two digraphs $G$ and $H$ is a digraph with vertex set $$V(G\square H)=V(G)\times V(H)=\{(x, x')\colon\, x\in V(G), x'\in V(H)\}$$ and arc set $$A(G\square H)=\{(x,x')(y,y')\colon\, xy\in A(G), x'=y',~or~x=y, x'y'\in A(H)\}.$$ 
We define the {\em $k$th powers with respect to Cartesian product} as $$D^{\square k}=\underbrace{D\square D\square \cdots \square D}_{k \mbox{ times}}.$$

The {\em strong product} $G\boxtimes H$ of two digraphs $G$ and $H$ is a digraph with vertex set $$V(G\boxtimes H)=V(G)\times V(H)=\{(x,
x')\colon\, x\in V(G), x'\in V(H)\}$$ and arc set $$A(G\boxtimes H)=A(G\square H)\cup \{(x,x')(y,y')\mid xy\in A(G),~x'y'\in A(H)\}.$$ 

The {\em lexicographic product} $G\circ H$ of two digraphs $G$ and $H$ is a digraph with vertex set $$V(G\circ H)=V(G)\times V(H)=\{(x,
x')\colon\, x\in V(G), x'\in V(H)\}$$ and arc set $$A(G\circ H)=\{(x,x')(y,y')\colon\, xy\in A(G),~or~x=y~and~x'y'\in
A(H)\}.$$

\section{Directed Steiner tree packing problem}\label{sec:tree}

As shown in the following two former subsections, the complexity for the decision version of {\sc IDSTP} and {\sc ADSTP} on general digraphs, symmetric digraphs and Eulerian digraphs have been completely determined. The hardness of approximation of both {\sc ADSTP} and {\sc IDSTP} will be mentioned. In the last subsection, we will also introduce three related topics: directed tree connectivity, independent directed Steiner trees and arc-disjoint branchings.

\subsection{General digraphs}

We first introduce the well-known problem of {\sc Directed $k$-Linkage} which is formulated as follows: 
\vspace{2mm}

\noindent {\bf Directed $k$-Linkage}: 
For a fixed integer $k\ge 2$, given a digraph $D$ and a (terminal) sequence $((s_1,t_1),\dots ,(s_k,t_k))$ of distinct vertices of $D,$ decide whether $D$ has $k$ vertex-disjoint paths $P_1,\dots ,P_k$, where $P_i$ starts at $s_i$ and ends at $t_i$ for all $i\in [k].$

\vspace{2mm}
Fortune, Hopcroft and Wyllie proved the following important theorem on 2-linkage problem.

\begin{thm}\label{thm2-14}\cite{Fortune-Hopcroft-Wyl}
The {\sc Directed 2-Linkage} is NP-complete.
\end{thm}

Let $D$ be a digraph and $S \subseteq V(D)$ with $|S|=k$. It is
natural to consider the following problem: what is the complexity of
deciding whether $\kappa_{S, r}(D)\ge \ell$~(resp. $\lambda_{S, r}(D)\ge
\ell)$? where $r\in S$ is a root. If $k=2$, say $S=\{r, x\}$, then
the problem of deciding whether $\kappa_{S, r}(D)\ge
\ell$~(resp.  $\lambda_{S, r}(D)\ge \ell)$ is equivalent to deciding whether
$\kappa(r, x)\ge \ell$~(resp.  $\lambda(r, x)\ge \ell)$, and so is
polynomial-time solvable (see \cite{Bang-Jensen-Gutin}), where
$\kappa(r, x)$~(resp. $\lambda(r, x))$ is the local vertex-strong
(resp. arc-strong) connectivity from $r$ to $x$. If $\ell=1$, then the above problem is also polynomial-time solvable by the well-known fact that every strong digraph has an out- and in-branching rooted at any vertex, and these branchings can be found in polynomial-time. Hence, it remains to consider the case that $k\geq 3,
\ell\geq 2$. Using Theorem~\ref{thm2-14}, Cheriyan and Salavatipour proved the NP-hardness of the case $k= 3, \ell= 2$ for both $\kappa_{S, r}(D)$ and $\lambda_{S, r}(D)$.

\begin{thm}\label{thm2-1}\cite{Cheriyan-Salavatipour}
Let $D$ be a digraph and $S \subseteq V(D)$ with $|S|=3$. The
problem of deciding whether $\kappa_{S, r}(D)\ge 2$ is NP-hard,
where $r\in S$.
\end{thm}

\begin{thm}\label{thm2-2}\cite{Cheriyan-Salavatipour}
Let $D$ be a digraph and $S \subseteq V(D)$ with $|S|=3$. The
problem of deciding whether $\lambda_{S, r}(D)\ge 2$ is NP-hard,
where $r\in S$.
\end{thm}

Sun and Yeo extended the above results to the following.

\begin{thm}\label{thm2-3}\cite{Sun-Yeo}
Let $k\geq 3$ and $\ell\geq 2$ be fixed integers (considered as constants). Let $D$ be a
digraph and $S \subseteq V(D)$ with $|S|=k$ and $r \in S$. Both the following problems are NP-complete.
\begin{itemize}
\item Is $\kappa_{S, r}(D)\ge \ell$?
\item Is $\lambda_{S, r}(D)\ge \ell$?
\end{itemize}
\end{thm}

The above results imply the entries in Tables~2-1 and~2-2.

\begin{center}
\begin{tabular}{|c||c|c|c|} \hline
\multicolumn{4}{|c|}{Table 2-1: General digraphs} \\ \hline
$\lambda_{S,r}(D) \geq \ell$? & $k=3$                                    & $k \geq 4$     &  $k$ part \\
$ |S|=k$                     &                                          & constant       &  of input \\ \hline  \hline
$\ell =2$                    & NP-complete \cite{Cheriyan-Salavatipour} & NP-complete\cite{Sun-Yeo}    &  NP-complete\cite{Sun-Yeo}  \\ \hline
$\ell \geq 3$ constant       & NP-complete \cite{Sun-Yeo}                             & NP-complete \cite{Sun-Yeo}   &  NP-complete\cite{Sun-Yeo}  \\ \hline
$\ell$ part of input         & NP-complete \cite{Sun-Yeo}                             & NP-complete\cite{Sun-Yeo}    &  NP-complete\cite{Sun-Yeo}  \\ \hline
\end{tabular}
\end{center}

\begin{center}
\begin{tabular}{|c||c|c|c|} \hline
\multicolumn{4}{|c|}{Table 2-2: General digraphs} \\ \hline
$\kappa_{S,r}(D) \geq \ell$? & $k=3$                                    & $k \geq 4$     &  $k$ part \\
$ |S|=k$                     &                                          & constant       &  of input \\ \hline  \hline
$\ell =2$                    & NP-complete \cite{Cheriyan-Salavatipour} & NP-complete\cite{Sun-Yeo}    &  NP-complete\cite{Sun-Yeo}  \\ \hline
$\ell \geq 3$ constant       & NP-complete \cite{Sun-Yeo}                             & NP-complete\cite{Sun-Yeo}    &  NP-complete\cite{Sun-Yeo}  \\ \hline
$\ell$ part of input         & NP-complete \cite{Sun-Yeo}                             & NP-complete\cite{Sun-Yeo}    &  NP-complete \cite{Sun-Yeo} \\ \hline
\end{tabular}
\end{center}

For general digraphs, Cheriyan and Salavatipour studied the hardness of approximation of both {\sc ADSTP} and {\sc IDSTP} as follows.

\begin{thm}\label{thm2-17}\cite{Cheriyan-Salavatipour}
Given an instance of ADSTP, it is NP-hard to approximate the solution within $O(m^{1/3-\epsilon})$ for any $\epsilon >0$.

\end{thm}

\begin{thm}\label{thm2-18}\cite{Cheriyan-Salavatipour}
Given an instance of IDSTP, it is NP-hard to approximate the solution within $O(n^{1/3-\epsilon})$ for any $\epsilon >0$.
\end{thm}

There are two more general problems:

\vspace{2mm}

\noindent {\bf GDE:}
The input consists of a digraph $D$, a capacity $c_e$ on each arc $e\in A(D)$, $\ell$ terminal sets $T_1, \dots, T_{\ell}$ and $\ell$ roots $r_1, \dots, r_{\ell}$ such that $r_i\in V(T_i)$ for each $i\in [\ell]$. The goal is to find a largest collection of directed Steiner trees, each rooted at an $r_i$ and containing all vertices of $T_i$ such that each arc $e$ is contained in at most $c_e$ directed trees.

\vspace{2mm}

\noindent {\bf GDV:}
The input consists of a digraph $D$, a capacity $c_v$ on each vertex $v\in V(D)$, $\ell$ terminal sets $T_1, \dots, T_{\ell}$ and $\ell$ roots $r_1, \dots, r_{\ell}$ such that $r_i\in V(T_i)$ for each $i\in [\ell]$. The goal is to find a largest collection of directed Steiner trees, each rooted at an $r_i$ and containing all vertices of $T_i$ such that each non-terminal vertex $v$ is contained in at most $c_v$ directed trees.

\vspace{2mm}

Clearly, when $\ell=1$ and $c_e=1$ (resp.  $c_v=1$), then {\sc GDE} (resp. {\sc GDV}) is exactly {\sc ADSTP} (resp. {\sc IDSTP}). It is worth mentioning that the hardness of approximation of both {\sc GDE} and {\sc GDV} has also been studied by researchers \cite{Cheriyan-Salavatipour, Guruswami2003}.

\subsection{Symmetric digraphs and Eulerian digraphs}

Now we turn our attention to symmetric digraphs. We first need the following theorem by Robertson and Seymour.

\begin{thm}\label{thm2-4} \cite{Robertson-Seymour}
Let $G$ be a graph and let $s_1,s_2,\ldots,s_r, t_1,t_2,\ldots,t_r$ be $2r$ distinct vertices in $G$.
We can in $O(|V(G)|^3)$ time decide if there exists an $(s_i,t_i)$-path, $P_i$, such that all $P_1,P_2,\ldots,P_r$ are vertex disjoint.
\end{thm}

Using Theorem~\ref{thm2-4}, Sun and Yeo proved the following result which implies the problem of {\sc Directed $k$-Linkage} can be solved in polynomial-time for symmetric digraphs, when $k$ is fixed.

\begin{cor}\label{thm2-5}\cite{Sun-Yeo}
Let $D$ be a symmetric digraph and let $s_1,s_2,\ldots,s_r$, $t_1,t_2,\ldots,t_r$ be vertices in $D$ (not necessarily distinct) and
let $S \subseteq V(D)$. We can in $O(|V(D)|^3)$ time decide if there for all $i=1,2,\ldots,r$ exists an $(s_i,t_i)$-path, $P_i$,
such that no internal vertex of any $P_i$ belongs to $S$ or to any path $P_j$ with $j \not=i$ (the end-vertices of $P_j$ can also not
be internal vertices of $P_i$).
\end{cor}

Furthermore, Sun and Yeo proved the following result on symmetric digraphs by Corollary~\ref{thm2-5}.

\begin{thm}\label{thm2-6}\cite{Sun-Yeo}
Let $k\geq 3$ and $\ell \geq 2$ be fixed integers and let $D$ be a symmetric digraph. Let $S \subseteq V(D)$ with $|S|=k$ and let $r$ be an arbitrary vertex in $S$. Let $A_0,A_1,A_2,\ldots , A_{\ell}$ be a partition of the arcs in $D[S]$.

We can in time $O(n^{\ell k - 2 \ell + 3} \cdot (2k-3)^{\ell (2k-3)})$
decide if there exist $\ell$ pairwise internally-disjoint $(S,r)$-trees, $T_1,T_2,\ldots ,T_{\ell}$, with $A(T_i) \cap A[S] = A_i$ for all $i=1,2,\ldots,\ell$ (note that $A_0$ are the arcs in $D[S]$ not used in any of the trees).
\end{thm}

The next corollary implies all the polynomial entries in Table~2-4.

\begin{cor}\label{thm2-7}\cite{Sun-Yeo}
Let $k\geq 3$ and $\ell \geq 2$ be fixed integers.
We can in polynomial time decide if $\kappa_{S,r}(D) \geq \ell$ for any symmetric
digraph, $D$, with $S \subseteq V(D)$, with $|S|=k$ and $ r\in S$.
\end{cor}

We now turn our attention to the NP-complete cases in Table~2-4. Chen, Li, Liu and Mao\cite{Chen-Li-Liu-Mao} introduced the following problem, which turned out to be NP-complete.

\vspace{2mm}

\noindent {\bf CLLM Problem:} Given a tripartite graph $G=(V, E)$ with a 3-partition $(A,B,C)$ such that $|A|=|B|=|C|=q$, decide whether there is a partition of $V$ into $q$ disjoint 3-sets $V_1, \dots, V_q$ such that for every $V_i= \{a_{i_1}, b_{i_2}, c_{i_3}\}$ $a_{i_1} \in A, b_{i_2} \in B, c_{i_3} \in C$ and $G[V_i]$ is connected.

\vspace{2mm}

\begin{lem}\label{thm2-8}\cite{Chen-Li-Liu-Mao} The  {\sc CLLM Problem} is NP-complete.
\end{lem}

Using Lemma~\ref{thm2-8}, Sun and Yeo showed the following. 

\begin{thm}\label{thm2-9}\cite{Sun-Yeo}
Let $k\geq 3$ be a fixed integer. The problem of deciding if a symmetric digraph $D$, with a $k$-subset $S$ of $V(D)$ with $r\in S$ satisfies $\kappa_{S, r}(D)\geq \ell$ ($\ell$ is part of the input), is NP-complete.
\end{thm}

The {\sc 2-Coloring Hypergraphs Problem} is defined as the following:

\vspace{2mm}

\noindent {\bf 2-Coloring Hypergraphs Problem:}
Given a hypergraph $H$ with vertex set $V(H)$ and edge set $E(H)$, determine if we can $2$-color the vertices $V(H)$ such that every hyperedge in $E(H)$ contains vertices of both colors. 

\vspace{2mm}

This problem is known to be NP-hard by Lov\'{a}sz.

\begin{thm}\label{2-Coloring-hyper}\cite{Lovasz}
The {\sc 2-Coloring Hypergraphs Problem} is NP-hard.
\end{thm}

By constructing a reduction from the {\sc 2-Coloring Hypergraphs Problem}, Sun and Yeo proved the following theorem.

\begin{thm}\label{thm2-10}\cite{Sun-Yeo}
Let $\ell \geq 2$ be a fixed integer. The problem of deciding if a symmetric digraph $D$, with an
$S \subseteq V(D)$ and $r\in S$ satisfies $\kappa_{S, r}(D)\geq \ell$ ($k=|S|$ is part of the input), is NP-complete.
\end{thm}

Theorem~\ref{thm2-9} together with Theorem~\ref{thm2-10} implies all the NP-completeness results in Table~2-4.


\begin{center}
		\begin{tabular}{|c|c|c|} \hline
			\multicolumn{3}{|c|}{Table 2-3: Symmetric digraphs} \\ \hline
			$\lambda_{S,r}(D) \geq \ell$? $ |S|=k$                                    & $k \geq 3$~constant     &  $k$ part of input\\
			                                                              \hline  \hline
			
			$\ell \geq 2$ constant                                     & Polynomial\cite{Sun-Yeo}    &  Polynomial\cite{Sun-Yeo}  \\ \hline
			$\ell$ part of input                                       & Polynomial\cite{Sun-Yeo}    &  Polynomial\cite{Sun-Yeo}  \\ \hline
		\end{tabular}
\end{center}


\begin{center}
		\begin{tabular}{|c|c|c|} \hline
			\multicolumn{3}{|c|}{Table 2-4: Symmetric digraphs} \\ \hline
			$\kappa_{S,r}(D) \geq \ell$? $ |S|=k$                                    & $k \geq 3$~constant     &  $k$ part of input\\
			                                                              \hline  \hline
			
			$\ell \geq 2$ constant                                     & Polynomial\cite{Sun-Yeo}    &  NP-complete\cite{Sun-Yeo}  \\ \hline
			$\ell$ part of input                                       & NP-complete\cite{Sun-Yeo}    &  NP-complete\cite{Sun-Yeo}  \\ \hline
		\end{tabular}
\end{center}

In a digraph $D$, the local arc-strong connectivity is the maximum number of arc-disjoint $(x,y)$-paths and is denoted by $\lambda_D(x,y)$, or just by $\lambda(x,y)$ if $D$ is clear from the context.

\begin{thm}\cite{JBJ-Frank-Jackson} \label{thm2-11}
Let $k \geq 1$ and let $D=(V,A)$ be a directed multigraph with a special vertex $r$.
Let $T' = \{ x \; | \; x \in V \setminus \{r\} \mbox{ and } d^-(x) < d^+(x) \}$.
If $\lambda(r,x) \geq k$ for every $x \in T'$, then there exists a family ${\cal F}$ of $k$ arc-disjoint out-trees rooted at $r$ so that every vertex
$x \in V$ belongs to at least $\min\{k,\lambda(r,x)\}$ members of ${\cal F}$.
\end{thm}

In the case when $D$ is Eulerian, then $d^+(x)=d^-(x)$ for all $x \in V(D)$ and therefore $T'=\emptyset$ in the above theorem.  Therefore the following
corollary holds.

\begin{cor} \label{thm2-12}\cite{Sun-Yeo}
Let $k \geq 1$ and let $D=(V,A)$ be an Eulerian digraph with a special vertex $r$.
Then there exists a family ${\cal F}$ of $k$ arc-disjoint out-trees rooted at $r$ so that every vertex
$x \in V$ belongs to at least $\min\{k,\lambda(r,x)\}$ members of ${\cal F}$.
\end{cor}

Using Corollary~\ref{thm2-12}, Sun and Yeo proved the following Theorem~\ref{thm2-13}. As one can determine $\lambda_D(r,s)$ in polynomial time for any $r$ and $s$ in $D$ we note that Theorem~\ref{thm2-13} implies that Table~2-5 (and therefore Table~2-3) holds.

\begin{thm}\label{thm2-13}\cite{Sun-Yeo}
If $D$ is an Eulerian digraph and $S \subseteq V(D)$ and $r \in S$, then $\lambda_{S,r}(D) \geq \ell$ if and only if $\lambda_D(r,s) \geq \ell$ for all $s \in S \setminus \{r\}$.
\end{thm}


\begin{center}
		\begin{tabular}{|c|c|c|} \hline
			\multicolumn{3}{|c|}{Table 2-5: Eulerian digraphs} \\ \hline
			$\lambda_{S,r}(D) \geq \ell$? $ |S|=k$                                    & $k \geq 3$~constant     &  $k$ part of input\\
			                                                              \hline  \hline
			
			$\ell \geq 2$ constant                                     & Polynomial\cite{Sun-Yeo}    &  Polynomial\cite{Sun-Yeo}  \\ \hline
			$\ell$ part of input                                       & Polynomial\cite{Sun-Yeo}    &  Polynomial\cite{Sun-Yeo}  \\ \hline
		\end{tabular}
\end{center}

Recently, Sun and Yeo proved that Theorem~\ref{thm2-14} holds even when $D$ is Eulerian.

\begin{thm}\label{thm2-15}\cite{Sun-Yeo}
The {\sc Directed 2-Linkage} restricted to Eulerian digraphs is NP-complete.
\end{thm}

Using Theorem~\ref{thm2-15}, we can prove the following result, which completes all the entries in Table~2-6.

\begin{thm}\label{thm2-16}\cite{Sun-Yeo}
Let $k\geq 3$ and $\ell\geq 2$ be fixed integers. Let $D$ be an Eulerian
digraph and $S \subseteq V(D)$ with $|S|=k$ and $r \in S$. Then deciding whether $\kappa_{S, r}(D)\ge \ell$ is NP-complete.
\end{thm}

\begin{center}
		\begin{tabular}{|c|c|c|} \hline
			\multicolumn{3}{|c|}{Table 2-6: Eulerian digraphs} \\ \hline
			$\kappa_{S,r}(D) \geq \ell$? $ |S|=k$                                    & $k \geq 3$~constant     &  $k$ part of input\\
			                                                              \hline  \hline
			
			$\ell \geq 2$ constant                                     & NP-complete\cite{Sun-Yeo}    &  NP-complete\cite{Sun-Yeo}  \\ \hline
			$\ell$ part of input                                       & NP-complete\cite{Sun-Yeo}    &  NP-complete\cite{Sun-Yeo}  \\ \hline
		\end{tabular}
\end{center}


It may be slightly surprising that for Eulerian digraphs the complexity of deciding if $\lambda_{S,r}(D) \geq \ell$ is always polynomial,
while the complexity of deciding if $\kappa_{S,r}(D) \geq \ell$ is always NP-complete.

It would also be interesting to determine the complexity for other classes of digraphs, like semicomplete digraphs. For example, one may consider the following question.

\begin{op}\label{op2-1}\cite{Sun-Yeo}
What is the complexity of deciding whether $\kappa_{S,r}(D)\ge \ell~$ (resp. $\lambda_{S,r}(D)\ge \ell)$
for integers $k\ge 3$ and $\ell\ge 2$, and a semicomplete digraph $D$?
\end{op}

In the argument for the question of deciding whether $\kappa_{S,r}(D)\ge \ell~$ for a symmetric digraph $D$ when both $k$ and $\ell$ are fixed, Sun and Yeo \cite{Sun-Yeo} used Corollary~\ref{thm2-5} where the $2k$ vertices $s_1, s_2, \ldots, s_k, t_1, t_2, \ldots, t_k$ are not necessarily distinct.
Normally in the $k$-linkage problem the initial and terminal vertices are considered distinct.  However if we allow them to be non-distinct and look
for internally-disjoint paths instead of vertex-disjoint paths, and the problem remains polynomial, then a similar approach to that of Theorem~\ref{thm2-6} and Corollary~\ref{thm2-7} can be used to show polynomiality of deciding whether $\kappa_{S,r}(D)\ge \ell~$.

\subsection{Related topics}


\subsubsection{Directed tree connectivity}

The following concept of directed tree connectivity is related to directed Steiner tree packing problem and is a natural
extension of tree connectivity \cite{Hager, Li-Mao5, Li-Mao-Sun} of undirected graphs to digraphs. The {\em generalized $k$-vertex-strong connectivity}
of $D$ is defined as $$\kappa_k(D)= \min \{\kappa_{S,r}(D)\mid S\subseteq V(D), |S|=k, r\in
S\}.$$ Similarly, the {\em generalized $k$-arc-strong connectivity} of $D$ is defined as $$\lambda_k(D)= \min \{\lambda_{S,r}(D)\mid S\subseteq V(D), |S|=k, r\in
S\}.$$ By definition, when $k=2$, $\kappa_2(D)=\kappa(D)$ and $\lambda_2(D)=\lambda(D)$. Hence, these two parameters could be seen as generalizations of classical vertex-strong connectivity and arc-strong connectivity of a digraph. The generalized $k$-vertex-strong connectivity and $k$-arc-strong connectivity are also called {\em directed tree connectivity}.

In \cite{Sun-Yeo}, Sun and Yeo proved some equalities and inequalities for directed tree connectivity.
They studied the relation between the directed tree connectivity and classical connectivity of digraphs by
showing that $\kappa_k(D)\leq \kappa(D)$ (when $n\geq k+\kappa(D)$) and $\lambda_k(D)=
\lambda(D)$. Furthermore, the upper bound for $\kappa_k(D)$ is sharp. Let $D$
be a strong digraph of order $n$. For $2\leq k\leq n$, they proved
that $1\leq \kappa_k(D)\leq n-1$ and $1\leq \lambda_k(D)\leq n-1$. All bounds are sharp, and they also characterized
those digraphs $D$ for which $\kappa_k(D)$~(respectively, $\lambda_k(D)$) attains
the upper bound. In the same paper, sharp Nordhaus-Gaddum type bounds for
$\lambda_k(D)$ were also given; moreover, extremal digraphs for the lower bounds were characterized.


A digraph $D=(V(D), A(D))$ is called {\em minimally generalized $(k,
\ell)$-vertex (resp. arc)-strongly connected} if $\kappa_k(D)\geq \ell$ (resp. $\lambda_k(D)\geq \ell$) but for any arc $e\in A(D)$, $\kappa_k(D-e)\leq \ell-1$ (resp. $\lambda_k(D-e)\leq \ell-1$), where $2\leq k\leq n, 1\leq \ell \leq n-1$. In \cite{Sun }, Sun studied the minimally generalized $(k, \ell)$-vertex-strongly connected digraphs and minimally generalized $(k, \ell)$-arc-strongly connected digraphs.




\subsubsection{Independent directed Steiner trees and independent branchings}

Let $D$ be a digraph with $r\in S\subseteq V(D)$.
In an $(S,r)$-tree $T$, we use $T(r,s)$ to denote the unique path in $T$ from the root $r$ to $s$, where $s\in S$. Two $(S,r)$-trees $T_1$ and $T_2$ are said {\em arc-independent}  if the paths $T_1(r,s)$ and $T_2(r,s)$ are arc-disjoint for every terminal $s\in S$. Two $(S,r)$-trees $T_1$ and $T_2$ are said {\em internally-independent}  if the paths $T_1(r,s)$ and $T_2(r,s)$ are internally-disjoint for each $s\in S$. If an $(S,r)$-tree $T$ is a spanning subdigraph of $D$, then it is called an {\em $r$-branching}  of $D$. Similarly, we can define the concepts of {\em arc-independent $r$-branchings}  and  {\em internally-independent $r$-branchings}. There are some publications on these topics, such as \cite{Berczi-Kovacs, Georgiadis-Tarjan, Huck-1995, Huck-1999a, Huck-1999b, Whitty}.

\subsubsection{Arc-disjoint in- and out-branchings rooted at the same vertex}


Recall that it is NP-complete to decide whether a digraph $D$ has a pair of arc-disjoint out-branching and in-branching rooted at $r,$ which was proved by Thomassen (for example, see \cite{Bang-Jensen}). Following \cite{Bang-Jensen-Huang2014} we will call such a pair a {\em good pair rooted at $r$}. Note that a good pair forms a strong spanning subdigraph of $D$ and thus if $D$ has a good pair, then $D$ is strong. The problem of the existence of a good pair was studied for tournaments and their generalizations, and characterizations (with proofs
implying polynomial-time algorithms for finding such a pair) were obtained in \cite{Bang-Jensen} for tournaments,  \cite{Bang-Jensen-Huang1995}
for quasi-transitive digraphs and \cite{Bang-Jensen-Huang2014} for locally semicomplete digraphs. Also, Bang-Jensen and Huang \cite{Bang-Jensen-Huang1995}  showed that if $r$ is adjacent to every vertex of $D$ (apart from itself) then $D$ has a good pair rooted at $r$. Gutin and Sun \cite{Gutin-Sun} studied the existence of a good pair in compositions of digraphs. They obtained the following result: every strong digraph composition $Q$ in which $n_i\ge 2$ for every $1\leq i\leq t$, has a good pair at every vertex of $Q.$ The condition of $n_i\ge 2$ in this result cannot be relaxed. They also characterizeed semicomplete compositions with a good pair, which generalizes the corresponding characterization by Bang-Jensen and Huang \cite{Bang-Jensen-Huang1995} for quasi-transitive digraphs. As a result, we can decide in polynomial time whether a given semicomplete composition has a good pair rooted at a given vertex. The readers can see \cite{Sun2} for a recent survey on the topic of digraph compositions.

\section{Directed Steiner path packing problem}\label{sec:steinerpath}

In this section, we will introduce the complexity for {\sc IDSPP} on Eulerian digraphs and symmetric digraphs, and the complexity for {\sc ADSPP} on general digraphs, when both $k$ and $\ell$ are fixed.

\subsection{Results for IDSPP}

By Theorem~\ref{thm2-15}, Sun and Zhang proved the NP-completeness of deciding whether $\kappa^p_{S,r}(D)\geq \ell$ for Eulerian digraphs (and therefore for general digraphs), which implies all entries of Table~3-1.

\begin{thm}\label{thm7-1-0}\cite{Sun5}
Let $k\geq 3, \ell \geq 1$ be fixed integers. For any Eulerian digraph $D$ and $S
\subseteq V(D)$ with $|S|=k$ and $r\in S$, the problem of deciding whether $\kappa^p_{S,r}(D) \geq \ell$ is NP-complete.
\end{thm}

	\begin{center}
		\begin{tabular}{|c|c|c|} \hline
			\multicolumn{3}{|c|}{Table 3-1: Eulerian digraphs (also for general digraphs)} \\ \hline
            $\kappa^p_{S,r}(D) \geq \ell$?~$|S|=k$                                    & $k \geq 3$~constant    &  $k$ part of input \\\hline  \hline
			
			$\ell \geq 1$ constant                                     & NP-complete\cite{Sun5}    &  NP-complete\cite{Sun5}  \\ \hline
			$\ell$ part of input                                       & NP-complete\cite{Sun5}   &  NP-complete\cite{Sun5}  \\ \hline
		\end{tabular}
	\end{center}

However, when we consider the class of symmetric digraphs, the problem becomes polynomial-time solvable according to the following theorem which can be deduced by Corollary~\ref{thm2-5}.

\begin{thm}\label{thm7-3-0}\cite{Sun5}
Let $k\geq 3$ and $\ell \geq 2$ be fixed integers. We can in polynomial time decide if $\kappa^p_{S,r}(D) \geq \ell$ for any symmetric digraph $D$ with $S \subseteq V(D)$, with $|S|=k$ and $r\in S$.
\end{thm}

Theorem~\ref{thm7-3-0} implies the polynomial entry of Table~3-2.

	\begin{center}
		\begin{tabular}{|c|c|c|} \hline
			\multicolumn{3}{|c|}{Table 3-2: Symmetric digraphs} \\ \hline
           $\kappa^p_{S,r}(D) \geq \ell$?~$|S|=k$                                    & $k \geq 3$~constant    &  $k$ part of input \\\hline  \hline
			
			$\ell \geq 2$ constant                                     & Polynomial\cite{Sun5}    &    \\ \hline
			$\ell$ part of input                                       &     &    \\ \hline
		\end{tabular}
	\end{center}

\begin{op}\label{op7-1-0}
Complete the blanks in Table~3-2.
\end{op}

It would also be interesting to study the complexity of {\sc IDSPP} on other digraph classes, such as semicomplete digraphs.

\begin{op}\label{op7-1}
Let $k\geq 3$ and $\ell \geq 1$ be fixed integers. What is the complexity of deciding whether $\kappa^p_{S,r}(D) \geq \ell$ for a semicomplete digraph $D$? where $r\in S \subseteq V(D)$ and $|S|=k$.
\end{op}

\subsection{Results for ADSPP}

The {\sc Directed Weak $k$-Linkage} problem is formulated as follows: 

\vspace{2mm}

\noindent {\bf Directed Weak $k$-Linkage}: For a fixed integer $k\ge 2$, given a digraph $D$ and a (terminal)
sequence $((s_1,t_1),\dots ,(s_k,t_k))$ of distinct vertices of $D,$ decide whether $D$ has $k$ arc-disjoint paths $P_1,\dots ,P_k$, where $P_i$ starts at $s_i$ and ends at $t_i$ for all $i\in [k].$ 

\vspace{2mm}
Fortune, Hopcroft and Wyllie proved the following important theorem on {\sc Directed Weak 2-Linkage}.

\begin{thm}\label{thm3-8}\cite{Fortune-Hopcroft-Wyl}
The {\sc Directed Weak 2-Linkage} is NP-complete.
\end{thm}

By Theorems~\ref{thm7-1-0} and~\ref{thm3-8}, Sun and Zhang proved the NP-completeness of deciding whether $\lambda^p_{S,r}(D)\geq \ell$ for general digraphs, which implies the entries of Table~3-3.

\begin{thm}\label{thm7-2}\cite{Sun5}
Let $k\geq 3, \ell \geq 1$ be fixed integers. For any digraph $D$ and $S
\subseteq V(D)$ with $|S|=k$ and $r\in S$, the problem of deciding whether $\lambda^p_{S,r}(D) \geq \ell$ is NP-complete.
\end{thm}

	\begin{center}
		\begin{tabular}{|c|c|c|} \hline
			\multicolumn{3}{|c|}{Table 3-3: General digraphs} \\ \hline
           $\lambda^p_{S,r}(D) \geq \ell$?~$|S|=k$                                    & $k \geq 3$~constant    &  $k$ part of input \\\hline  \hline
			
			$\ell \geq 1$ constant                                     & NP-complete\cite{Sun5}    &  NP-complete\cite{Sun5}  \\ \hline
			$\ell$ part of input                                       & NP-complete\cite{Sun5}    &  NP-complete\cite{Sun5}  \\ \hline
		\end{tabular}
	\end{center}  

It would also be interesting to study the complexity of {\sc ADSPP} on some digraph classes, such as semicomplete digraphs, Eulerian digraphs and symmetric digraphs.

\begin{op}\label{op7-2}
Let $k\geq 3$ and $\ell \geq 1$ be fixed integers. What is the complexity of deciding whether $\lambda^p_{S,r}(D) \geq \ell$ for a semicomplete digraph (Eulerian digraph, or symmetric digraph) $D$? where $r\in S \subseteq V(D)$ and $|S|=k$.
\end{op}

\subsection{A related topic: directed path connectivity}

Sun and Zhang~\cite{Sun5} introduced the following concept of directed path connectivity which is
related to directed Steiner path packing problem and is a natural
extension of path connectivity on undirected graphs (see \cite{Li-Mao5} for the introduction of path connectivity) to directed
graphs.
The {\em  directed path $k$-connectivity}
of $D$ is defined as
$$\kappa^p_k(D)= \min \{\kappa^p_{S,r}(D)\mid S\subseteq V(D), |S|=k, r\in
S\}.$$ Similarly, the {\em directed path $k$-arc-connectivity}
of $D$ is defined as
$$\lambda^p_k(D)= \min \{\lambda^p_{S,r}(D)\mid S\subseteq V(D), |S|=k, r\in
S\}.$$ By definition, when $k=2$, $\kappa^p_2(D)=\kappa(D)$ and
$\lambda^p_2(D)=\lambda(D)$. Hence, these two parameters could be seen
as generalizations of vertex-strong connectivity $\kappa(D)$ and arc-strong
connectivity $\lambda(D)$ of a digraph $D$. The directed path $k$-connectivity and the directed path $k$-arc-connectivity are also called {\em directed path connectivity}. 

In~\cite{Sun5}, Sun and Zhang showed that the values $\lambda^p_k(D)$ is decreasing over $k$, but the values $\kappa^p_k(D)$ are neither increasing, nor decreasing over $k$. They gave two upper bounds for $\kappa^p_k(D)$ in terms of $\kappa(D)$ and $\lambda(D)$. They also gave sharp lower and upper bounds for the Nordhaus-Gaddum type relations of the parameter $\lambda^p_k(D)$.

\section{Directed pendant Steiner tree packing problem}

\subsection{Results for IDPSTP}

In this section, we will introduce the complexity for {\sc IDPSTP} on Eulerian digraphs and symmetric digraphs.

By Theorem~\ref{thm2-15}, Yu and Sun~\cite{Yu-Sun2} showed that the problem of determining whether $\tau_{S,r}(D)\geq \ell$ on Eulerian digraphs (and therefore general digraphs) is NP-complete. This result implies the contents of Table~4-1.

\begin{thm}\label{thm7-1}\cite{Yu-Sun2}
Let $D$ be an Eulerian digraph and $S\subseteq V(D)$ with $|S|=k$ and $r\in S$. Let $k\geq 3$ and $\ell\geq 2$ be fixed integers. The problem of deciding if $\tau_{S,r}(D)\geq \ell$ is NP-complete.
\end{thm}

\begin{center}
		\begin{tabular}{|c|c|c|} \hline
			\multicolumn{3}{|c|}{Table 4-1: Eulerian digraphs ~(also for general digraphs)} \\ \hline
			$\tau_{S,r}(D) \geq \ell$? $ |S|=k$                                    & $k \geq 3$~constant     &  $k$ part of input\\
			                                                              \hline  \hline
			
			$\ell \geq 2$ constant                                     & NP-complete\cite{Yu-Sun2}    &  NP-complete\cite{Yu-Sun2}  \\ \hline
			$\ell$ part of input                                       & NP-complete\cite{Yu-Sun2}    &  NP-complete\cite{Yu-Sun2}  \\ \hline
		\end{tabular}
\end{center}


We now turn our attention to symmetric digraphs. Recall that a connected symmetric digraph is Eulerian. However, unlike Eulerian digraphs, the problem of deciding whether $\tau_{S,r}(D)\geq \ell$ on symmetric digraphs is polynomial-time solvable, when $|S|=k\geq 3$ and $\ell\geq 2$ are both fixed integers. By Corollary~\ref{thm2-5}, Yu and Sun~\cite{Yu-Sun2} gave the following theorem which shows the polynomiality for $\tau_{S,r}(D)$ on symmetric digraphs in Table~4-2.

\begin{thm}\label{thm7-2}\cite{Yu-Sun2}
    Let $D$ be a symmetric digraph and $S\subseteq V(D)$ with $|S|=k$ and $r\in S$. Let $k\geq 3$ and $\ell\geq 2$ be fixed integers. The problem of deciding if $\tau_{S,r}(D)\geq \ell$ can be solved in polynomial-time.
\end{thm}

We consider the NP-complete cases in Table~4-2 now. By Lemma~\ref{thm2-8}, Yu and Sun obtained the following result on symmetric digraphs. 


\begin{thm}\label{thm7-3}\cite{Yu-Sun2}
    Let $D$ be a symmetric digraph and $S\subseteq V(D)$ with $|S|=k$ and $r\in S$. Let $k\geq 3$ be a fixed integer. The problem of deciding if $\tau_{S,r}(D)\geq \ell$ ($\ell\geq 2$ is part of the input) is NP-complete.
\end{thm}

In order to prove the case that $k$ is part of input and $\ell\geq 2$ is a constant, Yu and Sun used Theorem~\ref{2-Coloring-hyper} and proved the following result, which (together with Theorem~\ref{thm7-3}) implies the NP-complete cases in Table~4-2.

\begin{thm}\label{thm7-4}\cite{Yu-Sun2}
    Let $D$ be a symmetric digraph and $S\subseteq V(D)$ with $|S|=k$ ($k$ is part of the input) and $r\in S$. Let $\ell\geq 2$ be a fixed integer. The problem of deciding if $\tau_{S,r}(D)\geq \ell$ is NP-complete.
\end{thm}

\begin{center}
		\begin{tabular}{|c|c|c|} \hline
			\multicolumn{3}{|c|}{Table 4-2: Symmetric digraphs} \\ \hline
			$\tau_{S,r}(D) \geq \ell$? $ |S|=k$                                    & $k \geq 3$~constant     &  $k$ part of input\\
			                                                              \hline  \hline
			
			$\ell \geq 2$ constant                                     & Polynomial \cite{Yu-Sun2}    &  NP-complete\cite{Yu-Sun2}  \\ \hline
			$\ell$ part of input                                       & NP-complete\cite{Yu-Sun2}    &  NP-complete\cite{Yu-Sun2}  \\ \hline
		\end{tabular}
\end{center}


For general digraphs, Yu and Sun studied the hardness of approximation of both {\sc IDPSTP} as follows.

\begin{thm}\label{thm7-5}\cite{Yu-Sun2}
Given an instance of {\sc IDPSTP}, it is NP-hard to approximate the solution within $O(n^{1/3-\epsilon})$ for any $\epsilon> 0$.
\end{thm}

\subsection{A related topic: directed pendant-tree connectivity}

Yu and Sun~\cite{Yu-Sun2} defined the directed pendant-tree $k$-connectivity of a digraph $D$ as
\begin{align*}
    \tau_{k}(D)=\min\{\tau_{S,r}(D)\mid S\subseteq V(D),|S|=k,r\in S\},
\end{align*}
where $\tau_{S,r}(D)$ denotes the maximum number of pairwise internally-disjoint pendant $(S,r)$-trees in $D$. By definition, it is clear that
\begin{align*}
    \begin{cases}
        \tau_{2}(D) = \kappa_{2}(D)=\kappa(D), \\
        \tau_{k}(D) \leq \kappa_{k}(D), 3\leq k\leq n.
    \end{cases}
\end{align*}
Hence, $\tau_{k}(D)$ is another type of generalization of the classical vertex-strong connectivity $\kappa(D)$ of digraphs. It is also worth noting that $\tau_k(D)$ is a directed version of the pendant-tree connectivity \cite{Hager}.

In \cite{Yu-Sun2}, Yu and Sun studied sharp bounds and values for the parameter $\tau_{k}(D)$. In particular, they gave some sharp bounds of $\tau_{k}(D)$ for a general digraph $D$. The first result is about a sharp lower bound and a sharp upper bound for $\tau_k(D)$ in terms of $n$ and $k$. The second result concerns a sharp upper bound for $\tau_{k}(D)$ in terms of arc-cuts of $D$. Finally, they studied Nordhaus-Gaddum type relations of $\tau_{k}(D)$ and obtained several sharp bounds of this type.

In \cite{Yu-Sun1}, Yu and Sun studied the directed pendant-tree 3-connectivity of symmetric compositions and derived a sharp lower bound for $\tau_{3}(Q)$, where $Q$ is a symmetric composition. Moreover, they proposed a polynomial-time algorithm for finding internally-disjoint pendant $(S,r)$-trees which attain the lower bound.

In \cite{Yu-Sun3}, Yu and Sun studied the directed pendant-tree 3-connectivity of Cartesian product digraphs $D\square H$, and derived a sharp lower bound for $\tau_{3}(D\square H)$, where $D$ and $H$ are both strong digraphs. Specifically, they proved the $\tau_{3}(D\square H)\geq \tau_{3}(D)+\tau_{3}(H)$. Moreover, they proposed a polynomial-time algorithm for finding internally-disjoint pendant $(S,r)$-trees which attain this lower bound.

\section{Strong subgraph packing problem}\label{sec:strongsubgraph}

As shown in the following two former subsections, we will introduce the complexity for the decision version of {\sc ISSP} or/and {\sc ASSP} on general digraphs, semicomplete digraphs, symmetric digraphs and Eulerian digraphs. Inapproximability results on {\sc  ISSP} and {\sc ASSP} will also be mentioned in this section. In the final subsection, we will introduce two related topics: the first one concerns Kriesell's conjecture and its extension in strong subgraph packing problem, and another one is about strong subgraph connectivity on digraphs, which is related to strong subgraph packing problem. 

\subsection{General digraphs}

For a fixed integer $k\ge 2$, it is easy to decide whether $\kappa_S(D)\ge 1$ for a digraph $D$ with $|S|=k$: it holds if and only if $D$ is strong. Unfortunately, deciding
whether $\kappa_S(D)\ge 2$ is already NP-complete for $S \subseteq V(D)$ with $|S|=k$, where $k\ge 2$ is a fixed integer.

By using the reduction from the {\sc Directed $k$-Linkage} problem (Theorem~\ref{thm2-14}), we can prove the following intractability result.

\begin{thm}\label{thm3-1}\cite{Sun-Gutin-Yeo-Zhang} Let $k\ge 2$ and $\ell\ge 2$ be fixed integers. Let $D$ be a digraph and $S \subseteq V(D)$ with $|S|=k$. The
problem of deciding whether $\kappa_S(D)\ge \ell$ is NP-complete.
\end{thm}

Yeo proved that it is an NP-complete problem to decide whether a 2-regular digraph has two arc-disjoint hamiltonian cycles (see,
e.g., Theorem 6.6 in \cite{Bang-Jensen-Yeo}). Thus, the problem of
deciding whether $\lambda_{V(D)}(D)\ge 2$ is NP-complete. Sun and Gutin \cite{Sun-Gutin2} extended this
result in Theorem~\ref{thm3-6} by Theorem~\ref{thm3-8}.

\begin{thm}\label{thm3-6}\cite{Sun-Gutin2}
Let $k\ge 2$ and $\ell\ge 2$ be fixed integers. Let $D$ be a digraph
and $S \subseteq V(D)$ with $|S|=k$. The problem of deciding whether
$\lambda_S(D)\ge \ell$ is NP-complete.
\end{thm}

Now Theorems~\ref{thm3-1} and \ref{thm3-6} imply the entries in Tables~5-1 and~5-2.

\begin{center}
\begin{tabular}{|c|c|c|} \hline
\multicolumn{3}{|c|}{Table 5-1: General digraphs} \\ \hline
$\kappa_{S}(D) \geq \ell$?                                     & $k \geq 2$     &  $k$ part \\
$ |S|=k$                                                              & constant       &  of input \\ \hline  \hline

$\ell \geq 2$ constant                                     & NP-complete \cite{Sun-Gutin-Yeo-Zhang}   &  NP-complete\cite{Sun-Gutin-Yeo-Zhang}  \\ \hline
$\ell$ part of input                                       & NP-complete\cite{Sun-Gutin-Yeo-Zhang}    &  NP-complete\cite{Sun-Gutin-Yeo-Zhang}  \\ \hline
\end{tabular}
\end{center}

\begin{center}
\begin{tabular}{|c|c|c|} \hline
\multicolumn{3}{|c|}{Table 5-2: General digraphs} \\ \hline
$\lambda_{S}(D) \geq \ell$?& $k \geq 2$     &  $k$ part \\
$ |S|=k$                                                              & constant       &  of input \\ \hline  \hline

$\ell \geq 2$ constant                                     & NP-complete \cite{Sun-Gutin2}   &  NP-complete\cite{Sun-Gutin2}  \\ \hline
$\ell$ part of input                                       & NP-complete\cite{Sun-Gutin2}    &  NP-complete\cite{Sun-Gutin2}  \\ \hline
\end{tabular}
\end{center}

\subsection{Semicomplete digraphs, symmetric digraphs and Eulerian digraphs}

Chudnovsky, Scott and Seymour \cite{Chud-Scott-Seymour} proved the following powerful result.

\begin{thm}\label{thm3-2}\cite{Chud-Scott-Seymour}
Let $k$ and $c$ be fixed positive integers. The  {\sc Directed
$k$-Linkage} problem on a digraph $D$ whose vertex set can be
partitioned into $c$ sets each inducing a semicomplete digraph and a
terminal sequence $((s_1,t_1),\dots ,(s_k,t_k))$ of  distinct
vertices of $D$ can be solved in polynomial time.
\end{thm}

The following nontrivial lemma can be deduced from Theorem~\ref{thm3-2}.

\begin{lem}\label{lem3-3}\cite{Sun-Gutin-Yeo-Zhang}
Let $k$ and $\ell$ be fixed positive integers. Let $D$ be a digraph
and let $X_1,X_2,\ldots,X_{\ell}$ be $\ell$ vertex disjoint subsets
of $V(D)$, such that $|X_i| \leq k$ for all $i\in [\ell]$. Let $X =
\cup_{i=1}^{\ell} X_i$ and assume that every vertex in $V(D)
\setminus X$ is adjacent to every other vertex in $D$. We can
in polynomial time decide if there exists vertex disjoint subsets
$Z_1,Z_2,\ldots,Z_{\ell}$ of $V(D)$, such that $X_i \subseteq Z_i$
and $D[Z_i]$ is strongly connected for each $i\in [\ell]$.
\end{lem}

Using Lemma \ref{lem3-3}, Sun, Gutin, Yeo and Zhang proved the following result for semicomplete digraphs.
\begin{thm}\label{thm3-4}\cite{Sun-Gutin-Yeo-Zhang}
Let $k\ge 2$ and $\ell\ge 2$ be fixed integers. Let $D$ be a digraph and $S \subseteq V(D)$ with $|S|=k$. The problem of deciding whether $\kappa_S(D)\ge \ell$ for a semicomplete digraph $D$ is polynomial-time solvable.
\end{thm}

For semicomplete digraphs, by Theorem~\ref{thm3-4}, we have the following entry of Table~5-3.

	\begin{center}
		\begin{tabular}{|c|c|c|} \hline
			\multicolumn{3}{|c|}{Table 5-3: Semicomplete digraphs} \\ \hline
			$\kappa_{S}(D) \geq \ell$? $ |S|=k$                                    & $k \geq 2$ constant    &  $k$ part of input \\ \hline  \hline
			
			$\ell \geq 2$ constant                                     & Polynomial\cite{Sun-Gutin-Yeo-Zhang}    &    \\ \hline
			$\ell$ part of input                                       &     &    \\ \hline
		\end{tabular}
	\end{center}


\begin{op}\label{op3-3}
Complete the blanks in Table~5-3.
\end{op}

The {\sc Directed $k$-Linkage} problem is polynomial-time solvable for planar digraphs \cite{Schr} and digraphs of bounded directed
treewidth \cite{JRST}. However, it seems that we cannot use the approach in proving Theorem \ref{thm3-4} directly as the structure of
minimum-size strong subgraphs in these two classes of digraphs is more complicated than in semicomplete digraphs. Certainly, we cannot
exclude the possibility that decide if $\kappa_S(D)\geq \ell$ in planar digraphs and/or in digraphs of bounded
directed treewidth is NP-complete.

\begin{op}\label{op3-1}
What is the complexity of deciding whether $\kappa_{S}(D)\ge \ell~$ for (fixed) integers $k\ge 2$ and $\ell\ge 2$, and a planar digraph $D$, where $S\subseteq V(D)$ with $|S|=k$?
\end{op}

\begin{op}\label{op3-2}
What is the complexity of deciding whether $\kappa_{S}(D)\ge \ell~$ for (fixed) integers $k\ge 2$ and $\ell\ge 2$, and a digraph $D$ of bounded
directed treewidth, where $S\subseteq V(D)$ with $|S|=k$?
\end{op}

It would be interesting to identify large classes of digraphs for which the $\kappa_S(D)\ge \ell$ problem can be decided in polynomial time.

Now restricted to symmetric digraphs $D$, for any fixed integer
$k\geq 3$, by Lemma~\ref{thm2-8}, the problem of deciding whether
$\kappa_S(D)\geq \ell~(\ell \geq 1)$ is NP-complete for $S\subseteq
V(D)$ with $|S|=k$.

\begin{thm}\label{thm3-5}\cite{Sun-Gutin-Yeo-Zhang}
For any fixed integer $k\geq 3$, given a symmetric digraph $D$, a
$k$-subset $S$ of $V(D)$ and an integer $\ell~(\ell \geq 1)$,
deciding whether $\kappa_S(D)\geq \ell$, is NP-complete.
\end{thm}

The last theorem assumes that $k$ is fixed but $\ell$ is a part of
input. When both $k$ and $\ell$ are fixed, the problem of deciding
whether $\kappa_S(D) \geq \ell$ for a symmetric digraph $D$, is
polynomial-time solvable. We will start with the following technical
lemma.

\begin{lem}\label{lem3-6}\cite{Sun-Gutin-Yeo-Zhang}
Let $k,\ell \geq 2$ be fixed. Let $G$ be a graph and let $S
\subseteq V(G)$ be an independent set in $G$ with $|S|=k$. For $i\in
[\ell]$, let $D_i$ be any set of arcs with both end-vertices in $S$.
Let a forest $F_i$ in $G$ be called $(S,D_i)$-{\em acceptable} if
the digraph $\overleftrightarrow{F_i}+D_i$ is strong and contains
$S$. In polynomial time, we can decide whether there exists
edge-disjoint forests $F_1,F_2,\ldots,F_{\ell}$ such that $F_i$ is
$(S,D_i)$-acceptable for all $i\in [\ell]$ and $V(F_i) \cap V(F_j)
\subseteq S$ for all $1 \leq i < j \leq \ell$.
\end{lem}

Sun, Gutin, Yeo and Zhang proved the following result by Lemma \ref{lem3-6}:

\begin{thm}\label{thm3-7}\cite{Sun-Gutin-Yeo-Zhang}
Let $k, \ell \geq 2$ be fixed. For any symmetric digraph $D$ and $S
\subseteq V(D)$ with $|S|=k$, we can in polynomial time decide
whether $\kappa_S(D) \geq \ell$.
\end{thm}

Recall that it was proved in \cite{Sun-Gutin-Yeo-Zhang} that $\kappa_2(\overleftrightarrow{G})=\kappa(G)$, which means that $\kappa_2(\overleftrightarrow{G})$ can be computed in polynomial time. In fact, the argument also means that $\kappa_{\{x, y\}}(\overleftrightarrow{G})=\kappa_{\{x,y\}}(G)$, that is, the maximum number of disjoint $x-y$ paths in $G$, therefore can be computed in polynomial time. Then combining with Theorems~\ref{thm3-5}, \ref{thm3-7} and~\ref{thm3-11} (which is deduced from Theorem~\ref{2-Coloring-hyper} by constructing a reduction from the {\sc 2-Coloring Hypergraphs Problem}), we can complete all the entries of Table~5-4.

\begin{thm}\label{thm3-11}\cite{Sun-Zhang} 
For any fixed integer $\ell\geq 2$, given a symmetric digraph $D$, a
$k$-subset $S$ of $V(D)$ and an integer $k~(k \geq 2)$,
deciding whether $\kappa_S(D)\geq \ell$, is NP-complete.
\end{thm}

\begin{center}
\begin{tabular}{|c||c|c|c|} \hline
\multicolumn{4}{|c|}{Table 5-4: Symmetric digraphs} \\ \hline
$\kappa_{S}(D) \geq \ell$? & $k=2$              & $k \geq 3$     &  $k$ part \\
$ |S|=k$                   &                    & constant       &  of input \\ \hline  \hline
$\ell \geq 2$ constant     & Polynomial \cite{Sun-Gutin-Yeo-Zhang}        & Polynomial  \cite{Sun-Gutin-Yeo-Zhang}   &  NP-complete \cite{Sun-Zhang} \\ \hline
$\ell$ part of input       & Polynomial \cite{Sun-Gutin-Yeo-Zhang}         & NP-complete \cite{Sun-Gutin-Yeo-Zhang}   &  NP-complete \cite{Sun-Gutin-Yeo-Zhang} \\ \hline
\end{tabular}
\end{center}

In the end of this subsection, we introduce results for Eulerian digraphs. Using Theorem~\ref{thm2-15}, Sun, Gutin and Zhang proved the following result for Eulerian digraphs which means all entries of Table~5-5.


\begin{thm}\label{thm3-12}\cite{Sun-Zhang}
Let $k, \ell \geq 2$ be fixed. For any Eulerian digraph $D$ and $S
\subseteq V(D)$ with $|S|=k$, deciding whether $\kappa_S(D) \geq \ell$ is NP-complete.
\end{thm}

\begin{center}
		\begin{tabular}{|c|c|c|} \hline
			\multicolumn{3}{|c|}{Table 5-5: Eulerian digraphs} \\ \hline
			$\kappa_{S}(D) \geq \ell$? $ |S|=k$                                     & $k \geq 2$ constant    &  $k$ part of input 
			 \\ \hline  \hline
			
			$\ell \geq 2$ constant                                     & NP-complete\cite{Sun-Zhang}    &  NP-complete\cite{Sun-Zhang}  \\ \hline
			$\ell$ part of input                                       & NP-complete\cite{Sun-Zhang}    &  NP-complete\cite{Sun-Zhang}  \\ \hline
		\end{tabular}
\end{center}

So far, there is still no result for $\lambda_S(D)$ on semicomplete digraphs, symmetric digraphs and Eulerian digraphs. Therefore, it would be interesting to consider the following question.

\begin{op}\label{op3-4}
What is the complexity of deciding whether $\lambda_{S}(D)\ge \ell$ for (fixed) integers $k\ge 3$ and $\ell\ge 2$, and a semicomplete digraph (symmetric digraph, Eulerian digraph) $D$?
\end{op}

\subsection{Inapproximability results on {\sc ISSP} and {\sc ASSP}}

In this subsection, we first introduce the following {\sc Set Cover Packing} problem:

\vspace{2mm}

\noindent {\bf Set Cover Packing}: 
Given a bipartite graph $G=(C\cup B, E)$, find a largest collection of pairwise disjoint set covers of $B$, where a {\em set cover} of $B$ is a subset $S\subseteq C$ such that each vertex of $B$ has a neighbor in $C$. 

\vspace{2mm}

Feige et al. \cite{Feige-Halldorsson-Kortsarz-Srinivasan} proved the following inapproximability result on the {\sc Set Cover Packing} problem.

\begin{thm}\label{thm3-9}\cite{Feige-Halldorsson-Kortsarz-Srinivasan}
Unless P=NP, there is no $o(\log{n})$-approximation algorithm for
{\sc Set Cover Packing}, where $n$ is the order of $G$.
\end{thm}

Sun, Gutin and Zhang got two inapproximability results on {\sc ISSP} and {\sc ASSP} by reductions
from the {\sc Set Cover Packing} problem.

\begin{thm}\label{thm3-10}\cite{Sun-Zhang}The following assertions hold:
\begin{description}
\item[(a) ] Unless P=NP, there is no $o(\log{n})$-approximation algorithm for {\sc ISSP}, even restricted to the case that $D$ is a symmetric
digraph and $S$ is independent in $D$, where $n$ is the order of $D$.
\item[(b) ] Unless P=NP, there is no $o(\log{n})$-approximation algorithm for {\sc ASSP}, even restricted to the case that $S$ is independent
in $D$, where $n$ is the order of $D$.
\end{description}
\end{thm}

\subsection{Related topics}

\subsubsection{Kriesell's conjecture and its extension in strong subgraph packing problem}

Let $G$ be a connected graph with $S\subseteq V(G)$. We say that a set of edges $C$ of $G$ an {\em $S$-Steiner-cut} if there are at least two components of $G-C$ which contain vertices of $S$. Similarly, let $D$ be a strong digraph and $S\subseteq V(D)$; we say that a set of arcs $C$ of $D$ an {\em $S$-strong subgraph-cut} if there are at least two strong components of $D-C$ which contain vertices of $S$.

Kriesell posed the following well-known conjecture which concerns an approximate min-max relation
between the size of an $S$-Steiner-cut and the number of edge-disjoint $S$-Steiner trees. 

\begin{conj}\label{kriesell}\cite{Kriesell}
Let $G$ be a graph and $S\subseteq V(G)$ with $|S|\geq 2$. If every
$S$-Steiner-cut in $G$ has size at least $2{\ell}$, then $G$
contains $\ell$ pairwise edge-disjoint $S$-Steiner trees.
\end{conj}

Lau \cite{Lau} proved that the conjecture holds if every
$S$-Steiner-cut in $G$ has size at least $26{\ell}$. West and Wu \cite{West-Wu} improved the bound significantly by showing that the conjecture still holds if $26{\ell}$ is replaced by $6.5{\ell}$. So far the best bound $5{\ell}+4$ was obtained by DeVos, McDonald and Pivotto as follows.  

\begin{thm}\label{thm5-501}\cite{DeVos-McDonald-Pivotto}
Let $G$ be a graph and $S\subseteq V(G)$ with $|S|\geq 2$. If every
$S$-Steiner-cut in $G$ has size at least $5{\ell}+4$, then $G$
contains $\ell$ pairwise edge-disjoint $S$-Steiner trees.
\end{thm}

Similar to Theorem~\ref{thm5-501}, it is natural to study an approximate min-max relation between the size of minimum $S$-strong subgraph-cut and the maximum number of arc-disjoint $S$-strong subgraphs in a digraph $D$. Here is an interesting problem which is analogous to
Conjecture~\ref{kriesell}.

\begin{op}\label{prob5-501}\cite{Sun-Zhang}
Let $D$ be a digraph and $S\subseteq V(D)$ with $|S|\geq 2$. Find a function
$f(\ell)$ such that the following holds: If every $S$-strong
subgraph-cut in $D$ has size at least $f(\ell)$, then $D$ contains
$\ell$ pairwise arc-disjoint $S$-strong subgraphs. 
\end{op}

Sun, Gutin and Zhang \cite{Sun-Zhang} proved that there is a linear function $f(\ell)$ for a strong symmetric digraph.

\subsubsection{Strong subgraph connectivity}

Related to strong subgraph packing problems, there is a type of connectivity on digraphs, called strong subgraph connectivity, including strong subgraph $k$-connectivity and strong subgraph $k$-arc-connectivity. The {\em strong subgraph $k$-connectivity} \cite{Sun-Gutin-Yeo-Zhang} is defined as
$$\kappa^s_k(D)=\min\{\kappa_S(D)\mid S\subseteq V(D), |S|=k\}.$$
Similarly, the {\em strong subgraph $k$-arc-connectivity} \cite{Sun-Gutin2} is defined as
$$\lambda^s_k(D)=\min\{\lambda_S(D)\mid S\subseteq V(D), |S|=k\}.$$

The strong subgraph connectivity is not only a natural extension of tree connectivity \cite{Hager, Li-Mao5, Li-Mao-Sun} of undirected graphs to digraphs, but also could be seen as a generalization of classical connectivity of undirected graphs by the fact that $\kappa^s_2(\overleftrightarrow{G})=\kappa(G)$ \cite{Sun-Gutin-Yeo-Zhang} and
$\lambda^s_2(\overleftrightarrow{G})=\lambda(G)$ \cite{Sun-Gutin2}. For more information on this topic, the reader can see \cite{Dong-Gutin-Sun, Sun-Gutin2, Sun-Gutin-Yeo-Zhang, Sun-Jin, Sun-Zhang} or a recent survey \cite{Sun-Gutin}.

\section{Strong arc decomposition problem}


Recall that it is NP-complete to decide whether a digraph has a strong arc decomposition (Theorem~\ref{thm4-13}), therefore it is natural to consider this problem on some digraph classes. In this section, we will introduce results on digraph classes, including semicomplete digraphs, locally semicomplete digraphs, split digraphs, digraph compositions and digraph products.

\subsection{Semicomplete digraphs, locally semicomplete digraphs and split digraphs}

Note that every digraph with a strong arc decomposition is 2-arc-strong. Bang-Jensen and Yeo characterized semicomplete digraphs with a strong arc decomposition.

\begin{figure}[!h]
\begin{center}
\tikzstyle{vertexX}=[circle,draw, top color=gray!5, bottom color=gray!30, minimum size=16pt, scale=0.6, inner sep=0.5pt]
\tikzstyle{vertexY}=[circle,draw, top color=gray!5, bottom color=gray!30, minimum size=20pt, scale=0.7, inner sep=1.5pt]

\begin{tikzpicture}[scale=0.4]
 \node (v1) at (1.0,6.0) [vertexX] {$v_1$};
 \node (v2) at (7.0,6.0) [vertexX] {$v_2$};
 \node (v3) at (1.0,1.0) [vertexX] {$v_3$};
 \node (v4) at (7.0,1.0) [vertexX] {$v_4$};
\draw [->, line width=0.03cm] (v1) -- (v2);
\draw [->, line width=0.03cm] (v2) -- (v3);
\draw [->, line width=0.03cm] (v3) -- (v4);
\draw [->, line width=0.03cm] (v4) -- (v1);
\draw [->, line width=0.03cm] (v1) to [out=285, in=75] (v3);
\draw [->, line width=0.03cm] (v3) to [out=105, in=255] (v1);
\draw [->, line width=0.03cm] (v2) to [out=285, in=75] (v4);
\draw [->, line width=0.03cm] (v4) to [out=105, in=255] (v2);
\end{tikzpicture}
\caption{Digraph $S_4$}\label{S4fig}
\end{center}
\end{figure}
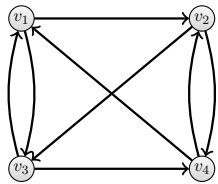

\begin{thm}\label{thm4-14}\cite{Bang-Jensen-Yeo}
A 2-arc-strong semicomplete digraph $D$ has a strong arc decomposition if and only if $D$ is not isomorphic to $S_4$, where $S_4$ is obtained from the complete digraph with four
vertices by deleting a cycle of length 4 (see Figure \ref{S4fig}). Furthermore, a strong arc decomposition of $D$ can be obtained in polynomial time when it exists.
\end{thm}

The following result extends Theorem \ref{thm4-14} to locally semicomplete digraphs.

\begin{thm}\label{thm4-15}\cite{Bang-Jensen-Huang2012}
A 2-arc-strong locally semicomplete digraph $D$ has a strong arc decomposition if and only if $D$ is not
the second power of an even cycle.
\end{thm}

Recently, Bang-Jensen and Wang studied strong arc decompositions of split digraphs and their main result in as follows:

\begin{thm}\label{thm4split}\cite{Bang-Jensen-Wang2025}
Let $D=(V_1,V_2;A)$ be a 2-arc-strong split digraph such that $V_1$ is an independent set and the subdigraph induced by $V_2$ is semicomplete. If every vertex of $V_1$ has both in- and out-degree at least three in $D$, then $D$ has a strong arc decomposition. Furthermore, every 3-arc-strong split digraph has a strong arc decomposition.
\end{thm}

\subsection{Digraph compositions}

Sun, Gutin and Ai \cite{Sun-Gutin-Ai} studied strong arc decompositions on digraph compositions and products.

The following theorem gives sufficient conditions for a digraph composition to have a strong arc decomposition.

\begin{thm}\label{thm4-1}\cite{Sun-Gutin-Ai}
Let $T$ be a digraph with vertices $u_1,\dots ,u_t$ ($t\ge 2$) and let $H_1,\dots ,H_t$ be digraphs. Let the vertex set of $H_i$ be $\{u_{i,j_i}\colon\, j_i\in [n_i]\}$ for every $i\in [t].$
Then $Q=T[H_1,\dots ,H_t]$ has a strong arc decomposition if at least one of the following conditions holds:
\begin{description}
\item[(a) ]  $T$ is a 2-arc-strong semicomplete digraph and $H_1,\dots , H_t$ are arbitrary digraphs, but $Q$ is not isomorphic to $S_4;$
\item[(b) ] $T$ has a Hamiltonian cycle and one of the following conditions holds:
\begin{itemize}
\item $t$ is even and $n_i\ge 2$ for every $i\in [t];$
\item $t$ is odd, $n_i\ge 2$ for every $i\in [t]$ and at least two distinct subdigraphs $H_i$ have arcs;
\item  $t$ is odd and $n_i\ge 3$ for every $i\in [t]$ apart from one $i$ for which $n_i\ge 2$.
\end{itemize}
\item[(c)] $T$ and all $H_i$ are strong digraphs with at least two vertices.
\end{description}
\end{thm}

They used Theorem~\ref{thm4-1} to prove the following characterization for semicomplete compositions $T[H_1,\dots ,H_t]$ when each $H_i$ has at least two vertices.

\begin{thm}\label{thm4-2}\cite{Sun-Gutin-Ai}
Let $Q=T[H_1,\dots ,H_t]$ be a strong semicomplete composition such that
each $H_i$ has at least two vertices for $i\in [t]$.
Then $Q$ has a strong arc decomposition if and only if $Q$ is not isomorphic to one of the following three digraphs:
$\overrightarrow{C}_3[\overline{K}_2,\overline{K}_2,\overline{K}_2]$, $\overrightarrow{C}_3[\overrightarrow{P_2},\overline{K}_2,\overline{K}_2]$,
$\overrightarrow{C}_3[\overline{K}_2,\overline{K}_2,\overline{K}_3].$
\end{thm}

Bang-Jensen, Gutin and Yeo solved an open problem in \cite{Sun-Gutin-Ai} by obtaining a characterization of all semicomplete compositions with a strong arc decomposition. 

\begin{thm}\label{thm4-3}\cite{Bang-Jensen-Gutin-Yeo}
Let $Q=T[H_1,\dots ,H_t]$ be a semicomplete composition. Then $Q$ has a strong arc decomposition if and only if $Q$ is 2-arc-strong  and is not isomorphic to one of the following four digraphs: $S_4$, $\overrightarrow{C}_3[\overline{K}_2,\overline{K}_2,\overline{K}_2]$, $\overrightarrow{C}_3[\overrightarrow{P_2},\overline{K}_2,\overline{K}_2]$,
$\overrightarrow{C}_3[\overline{K}_2,\overline{K}_2,\overline{K}_3].$
\end{thm}

Among their argument, the main technical result is the following theorem.

\begin{thm}\label{thm4-4}\cite{Bang-Jensen-Gutin-Yeo}
Let $Q=T[\overline{K}_{n_1}, \dots,\overline{K}_{n_t}]$ be an extended semicomplete digraph where $n_i\leq 2$ for $i\in [t]$. If $Q$ is 2-arc-strong, then $Q$ has a strong arc decomposition if and only if $Q$ is not isomorphic to one of the following four digraphs: $S_4$,
$\overrightarrow{C}_3[\overline{K}_2,\overline{K}_2,\overline{K}_2]$, $\overrightarrow{C}_3[\overrightarrow{P_2},\overline{K}_2,\overline{K}_2]$.
\end{thm}

The following theorem by Bang-Jensen and Huang gives a complete characterization of quasi-transitive digraphs and the decomposition below is called the {\em canonical decomposition} of a quasi-transitive digraph.

\begin{thm}\label{thm4-5}\cite{Bang-Jensen-Huang1995}
Let $D$ be a quasi-transitive digraph. The following assertions hold:
\begin{description}
\item[(a)] If $D$ is not strong, then there exists a transitive oriented
graph $T$ with vertices $\{u_i\mid i\in [t]\}$ and strong quasi-transitive digraphs $H_1, H_2, \dots, H_t$ such that $D = T[H_1, H_2, \dots, H_t]$, where $H_i$ is substituted for $u_i$, $i\in [t]$.
\item[(b)] If $D$ is strong, then there exists a strong semicomplete
digraph $S$ with vertices $\{v_j\mid j\in [s]\}$ and quasi-transitive digraphs $Q_1, Q_2, \dots, Q_s$ such that $Q_j$ is either a vertex or is non-strong and $D = S[Q_1, Q_2, \dots, Q_s]$, where $Q_j$ is substituted for $v_j$, $j\in [s]$.
\end{description}
\end{thm}

Theorems~\ref{thm4-3} and~\ref{thm4-5} imply a characterization of quasi-transitive digraphs with a strong arc decomposition (this solves another open question in \cite{Sun-Gutin-Ai}).

\begin{thm}\label{sss5}\cite{Bang-Jensen-Gutin-Yeo}
A quasi-transitive digraph $D$ has a strong arc decomposition if and only if $D$ is 2-arc-strong and is not isomorphic to one of the following four digraphs: $S_4$,
$\overrightarrow{C}_3[\overline{K}_2,\overline{K}_2,\overline{K}_2]$, $\overrightarrow{C}_3[\overrightarrow{P_2},\overline{K}_2,\overline{K}_2]$,
$\overrightarrow{C}_3[\overline{K}_2,\overline{K}_2,\overline{K}_3].$
\end{thm}

All proofs in \cite{Bang-Jensen-Gutin-Yeo} are constructive and can be turned into polynomial algorithms for finding strong arc decompositions. Thus, the problem of finding a strong arc decomposition in a semicomplete composition, which has one, admits a polynomial time algorithm.

By the definitions and Theorem~\ref{thm4-5}, strong semicomplete compositions generalize both strong semicomplete digraphs and strong quasi-transitive digraphs. However, they
do not generalize locally semicomplete digraphs and their generalizations in- and out-locally semicomplete digraphs. While there is a characterization of locally semicomplete digraphs having a strong arc decomposition \cite{Bang-Jensen-Huang2012}, no such a characterization is known for locally in-semicomplete digraphs and it would be interesting to obtain such a characterization or at least establish the complexity of deciding whether a locally in-semicomplete digraph has a strong arc decomposition.

\begin{op}\label{op4-3}\cite{Bang-Jensen-Gutin-Yeo}
Can we decide in polynomial time whether a given locally in-semicomplete digraph has a strong arc decomposition?
\end{op}

\begin{op}\label{op4-4}\cite{Bang-Jensen-Gutin-Yeo}
Characterize locally in-semicomplete digraphs with a strong arc decomposition.
\end{op}



Recall that a digraph $D$ has a strong arc decomposition if and only if $\lambda_{V(D)}(D)\geq 2$ (or, $\kappa_{V(D)}(D)\geq 2$). Therefore, the problem of {\sc ASSP} (or {\sc ISSP}) could be seen as an extension of the problem of strong arc decomposition. Sun, Gutin and Zhang gave two sufficient conditions for a digraph composition to have at least $n_0$ arc-disjoint $S$-strong subgraphs for any $S\subseteq V(Q)$ with $2\leq |S|\leq |V(Q)|$. Recall that $n_0=\min\{n_i\mid i\in [t]\}$.

\begin{thm}\label{thm4-16}\cite{Sun-Zhang}
The digraph composition $Q=T[H_1,\dots, H_t]$ has at least $n_0$ arc-disjoint $S$-strong subgraphs for any $S\subseteq V(Q)$ with $2\leq |S|\leq |V(Q)|$, if one of the following conditions holds:
\begin{description}
\item[(a) ] $D$ is a strong symmetric digraph;
\item[(b) ] $D$ is a strong semicomplete digraph and $Q\not\in
\mathcal{Q}_0$, where \\ $\mathcal{Q}_0=\{\overrightarrow{C}_3[\overline{K}_2,\overline{K}_2,\overline{K}_2], \overrightarrow{C}_3[\overrightarrow{P}_2,\overline{K}_2,\overline{K}_2], \overrightarrow{C}_3[\overline{K}_2,\overline{K}_2,\overline{K}_3]\}$.
\end{description}
Moreover, these strong subgraphs can be found within time
complexity $O(n^4)$, where $n$ is the order of $Q$.
\end{thm}

By Theorems~\ref{thm4-5} and~\ref{thm4-16}, we directly have:

\begin{cor}\label{thm4-12}\cite{Sun-Zhang}
Let $Q\not\in \mathcal{Q}_0$ be a strong quasi-transitive digraph. We can in polynomial time find at least $n_0$ arc-disjoint $S$-strong subgraphs in $Q$ for
any $S\subseteq V(Q)$ with $2\leq |S|\leq |V(Q)|$.
\end{cor}

\subsection{Digraph products}

Sun, Gutin and Ai obtained the following result on $G^{\square k}$ for any integer $k\geq 2$.

\begin{thm}\label{thm4-6}\cite{Sun-Gutin-Ai} Let $G$ be a strong digraph of order at least two which has a collection of arc-disjoint cycles covering all its
vertices and let $k\ge 2$ be an integer. The product digraph $D=G^{\square k}$ has a strong arc decomposition. Moreover, for any fixed integer $k$, such a strong arc decomposition can be found in polynomial time.
\end{thm}

They also obtained the following result on strong arc decomposition of strong product digraphs.

\begin{thm}\label{thm4-7}\cite{Sun-Gutin-Ai}
For any strong digraphs $G$ and $H$ with at least two vertices, the product digraph $D=G\boxtimes H$ has a strong arc decomposition. Moreover, such a decomposition can be found in polynomial time.
\end{thm}

By definition, the strong product digraph $G\boxtimes H$ is a spanning subdigraph of the lexicographic product digraph $G\circ H$, so the following result holds by Theorem \ref{thm4-7}: For any strong digraphs $G$ and $H$ with orders at least two, the product digraph $D=G\circ H$ has a strong arc decomposition. Moreover, such a decomposition can be found in polynomial time. In fact, we can get a more general result.

A digraph is {\em Hamiltonian decomposable} if it has a family of
Hamiltonian dicycles such that every arc of the digraph belongs to
exactly one of the dicycles. Ng \cite{Ng2} gives the most complete
result among digraph products.

\begin{thm}\label{thm4-8}\cite{Ng2}
If $G$ and $H$ are Hamiltonian decomposable digraphs, and $|V(G)|$
is odd, then $G\circ H$ is Hamiltonian decomposable.
\end{thm}

Theorem \ref{thm4-8} implies that if $G$ and $H$ are Hamiltonian
decomposable digraphs, and $|V(G)|$ is odd, then $G\circ H$ has a strong arc decomposition. It is
not hard to extend this result as follows: for any strong digraphs $G$ and $H$ of orders at least two, if $H$ contains $\ell \geq 1$ arc-disjoint strong spanning subdigraphs, then the product digraph $D=G\circ H$ can be decomposed into $\ell +1$ arc-disjoint strong spanning subdigraphs.

\begin{thm}\label{thm4-9}\cite{Trotter-Erdos}
The Cartesian product $\overrightarrow{C_p}\square
\overrightarrow{C_q}$ is Hamiltonian if and only if there are
non-negative integers $d_1, d_2$ for which $d_1+d_2= \gcd(p, q)\geq
2$ and $\gcd(p, d_1) = \gcd(q, d_2) = 1$.
\end{thm}

It was proved in \cite{Sun-Gutin-Ai} that $G\square H$ has a strong arc decomposition when $G\cong H$. However,
Theorem~\ref{thm4-9} implies this cannot be extended to the case that $G\not\cong H$, since it is not hard to show that the Cartesian product digraph of any two
cycles has a strong arc decomposition if and only if it has a pair of arc-disjoint Hamiltonian cycles. But
this could hold for the case that $G\not\cong H$ if we add additional conditions, since it was also showed in \cite{Sun-Gutin-Ai} that $G\square H$ has a strong arc decomposition when one of $G$ and $H$ has a strong arc decomposition. So the following open question is interesting:

\begin{op}\label{op4-1}\cite{Sun-Gutin-Ai}
For any two strong digraphs $G$ and $H$, neither of which have a strong arc decomposition, under what condition the product digraph $G\square H$ has a strong arc decomposition?
\end{op}

Furthermore, we may also consider the following more challenging question:

\begin{op}\label{op4-2}\cite{Sun-Gutin-Ai}
Under what condition the product digraph $G\square H~(G\boxtimes H)$ has more (than two) arc-disjoint strong spanning subdigraphs?
\end{op}

\subsection{Related topics: Kelly's conjecture and Bang-Jensen--Yeo conjecture}

The following is the famous Kelly's conjecture (see e.g. Conjecture~13.4.5 of \cite{Bang-Jensen-Gutin}):

\begin{conj}\label{Kelly}
Every regular tournament on $n=2k +1$ vertices has a decomposition into $k$ arc-disjoint Hamiltonian cycles.
\end{conj}

This conjecture was proved for large $n$ by K\"{u}hn and Osthus \cite{Kuhn-Osthus2013}. It is easy to see that a $k$-regular tournament is $k$-arc-strong. Bang-Jensen and Yeo posed the following conjecture which contains Conjecture~\ref{Kelly} as the special case when $n=2k+1$.

\begin{conj}\label{conj4-1-0}\cite{Bang-Jensen-Yeo}
Every $k$-arc-strong tournament has an arc-decomposition into $k$ arc-disjoint spanning strong subgraphs. Furthermore, for every natural number $k$, there exists a natural number $n_k$ such that every $k$-arc-strong semicomplete digraph with at least $n_k$ vertices has an arc-decomposition into $k$ arc-disjoint spanning strong subgraphs.
\end{conj}

Conjecture~\ref{conj4-1-0} is equivalent to the following: for a tournament $T$, if  $T$ is $k$-arc-strong, then $\kappa_{V(T)}(T)\geq k$. Bang-Jensen and Yeo  proved three results which support the conjecture, the first one is Theorem~\ref{thm4-14} and the remaining two are as follows:

\begin{thm}\label{thm4-10}\cite{Bang-Jensen-Yeo}
Every tournament which has a non-trivial cut (both sides containing at least two vertices) with precisely $k$ arcs in one direction contains $k$ arc-disjoint spanning strong subgraphs. 
\end{thm}

\begin{thm}\label{thm4-11}\cite{Bang-Jensen-Yeo}
Every $k$-arc-strong tournament with minimum semi-degree at least $37k$ contains $k$ arc-disjoint spanning strong subgraphs.
\end{thm}

Theorem~\ref{thm4-11} implies that if $T$ is a $74k$-arc-strong tournament with specified not necessarily distinct vertices $u_1, u_2, \dots, u_k, v_1, v_2, \dots, v_k$ then $T$ contains $2k$ arc-disjoint branchings $B^-_{u_1}, B^-_{u_2}, \dots, B^-_{u_k}$,  $B^+_{v_1}, B^+_{v_2}, \dots, B^+_{v_k}$ where $B^-_{u_i}$  is an in-branching rooted at the vertex $u_i$ and $B^+_{v_i}$ is an out-branching rooted at the vertex $v_i$, $i\in [k]$. This solves a conjecture of Bang-Jensen and Gutin \cite{Bang-Jensen-Gutin1998}.

Conjecture~\ref{conj4-1-0} does not hold for locally semicomplete digraphs as by Theorem~\ref{thm4-15} the second
power of an even cycle cannot be decomposed into two arc-disjoint Hamiltonian cycles. If $n$ is relatively
prime to both 2 and 3, then it is easy to see that $\overrightarrow{C}^3_n$
can be decomposed into three arc-disjoint
Hamiltonian cycles, where $\overrightarrow{C}^3_n$ is the 3rd power of $\overrightarrow{C}_n$. In fact, such a decomposition does not exist for any other $n$.

\begin{pro}\cite{Bang-Jensen-Huang2012}
If $n$ is not relatively prime to 2 or 3, then $\overrightarrow{C}^3_n$ cannot be decomposed into arc-disjoint Hamiltonian cycles.
\end{pro}

\section{Directed Steiner cycle packing problem}

In this section, we will introduce the complexity for {\sc IDSCP} on Eulerian digraphs and symmetric digraphs, and the complexity for {\sc ADSCP} on general digraphs, when both $k$ and $\ell$ are fixed. In the final part, we will introduce a related topic: directed cycle connectivity.

\subsection{Results for IDSCP}

Using Theorem~\ref{thm2-15}, Sun and Jin proved the following result on the NP-completeness of deciding whether $\kappa^c_S(D)\geq \ell$ for Eulerian digraphs (and therefore for general digraphs), which implies all entries of Table~7-1.

\begin{thm}\label{thm5-1}\cite{Sun4}
Let $k\geq 2, \ell \geq 1$ be fixed integers. For any Eulerian digraph $D$ and $S
\subseteq V(D)$ with $|S|=k$, deciding whether $\kappa^c_S(D) \geq \ell$ is NP-complete.
\end{thm}


\begin{center}
		\begin{tabular}{|c|c|c|} \hline
			\multicolumn{3}{|c|}{Table 7-1: Eulerian digraphs (also for general digraphs)} \\ \hline
			$\kappa^c_{S}(D) \geq \ell$? $ |S|=k$                                    & $k \geq 2$~constant     &  $k$ part of input\\
			                                                              \hline  \hline
			
			$\ell \geq 1$ constant                                     & NP-complete\cite{Sun4}    &  NP-complete\cite{Sun4}  \\ \hline
			$\ell$ part of input                                       & NP-complete\cite{Sun4}    &  NP-complete\cite{Sun4}  \\ \hline
		\end{tabular}
	\end{center}

When $D$ is an Eulerian digraph, we show in Theorem~\ref{thm5-1} that the problem of deciding whether $\kappa^c_S(D)\geq \ell$ with $|S|=k$ is NP-complete, where both $k\geq 2, \ell\geq 1$ are fixed integers. However, when we consider the class of symmetric digraphs, the problem becomes polynomial-time solvable, as shown in the following theorem which can be deduced from Corollary~\ref{thm2-5} and implies the polynomial entry of Table~7-2.  




\begin{thm}\label{thm5-3}\cite{Sun4}
Let $k\geq 2$ and $\ell \geq 1$ be fixed integers. We can in polynomial time decide if $\kappa^c_{S}(D) \geq \ell$ for any symmetric
digraph $D$ with $S \subseteq V(D)$, where $|S|=k$.
\end{thm}

\begin{center}
		\begin{tabular}{|c|c|c|} \hline
			\multicolumn{3}{|c|}{Table 7-2: Symmetric digraphs} \\ \hline
			$\kappa^c_{S}(D) \geq \ell$? $ |S|=k$                                     & $k \geq 2$~constant    &  $k$ part of input\\
			                                                               \hline  \hline
			
			$\ell \geq 1$ constant                                     & Polynomial\cite{Sun4}    &    \\ \hline
			$\ell$ part of input                                       &     &    \\ \hline
		\end{tabular}
	\end{center}

\begin{op}\label{op5-1-0}
Complete the blanks in Table~7-2.
\end{op}

We directly have the following result by Theorem~\ref{thm5-3}.

\begin{cor}\label{cor5-4}\cite{Sun4}
We can in polynomial time decide if a given symmetric digraph $D$ is $k$-cyclic, for a fixed integer $k$.
\end{cor}

It would also be interesting to study the complexity of {\sc IDSCP} on other digraph classes, such as semicomplete digraphs.

\begin{op}\label{op5-1}
Let $k\geq 2$ and $\ell \geq 1$ be fixed integers. What is the complexity of deciding whether $\kappa^c_{S}(D) \geq \ell$ for a semicomplete digraph $D$? where $S \subseteq V(D)$ and $|S|=k$.
\end{op}

\subsection{Results for ADSCP}

In this subsection, we turn our attention to the complexity for $\lambda^c_S(D)$ and obtain the following result on general digraphs by Theorems~\ref{thm2-15} and~\ref{thm5-1}, which implies all entries of Table~7-3.

\begin{thm}\label{thmd}\cite{Sun4}
Let $k\geq 2, \ell \geq 1$ be fixed integers. For a digraph $D$ and $S
\subseteq V(D)$ with $|S|=k$, deciding whether $\lambda^c_S(D) \geq \ell$ is NP-complete.
\end{thm}

	\begin{center}
		\begin{tabular}{|c|c|c|} \hline
			\multicolumn{3}{|c|}{Table 7-3: General digraphs} \\ \hline
			$\lambda^c_{S}(D) \geq \ell$? $ |S|=k$                                     & $k \geq 2$~constant   &  $k$ part of input\\
\hline  \hline
			
			$\ell \geq 1$ constant                                     & NP-complete\cite{Sun4}    &  NP-complete\cite{Sun4}  \\ \hline
			$\ell$ part of input                                       & NP-complete\cite{Sun4}    &  NP-complete\cite{Sun4}  \\ \hline
		\end{tabular}
	\end{center}    


Now we focus on digraph classes, including Eulerian digraphs, planar digraphs and symmetric digraphs. 
The following Flow Decomposition Theorem will be used in our argument.
\begin{thm}\label{07-0}[Flow Decomposition Theorem]\label{Flow Decomposition Theorem}\cite{Bang-Jensen-Gutin}
Every flow $x$ in $\mathcal{N}$ can be represented as the arc sum of some path and cycle flows $$f(P_1), f(P_2), \dots, f(P_\alpha), f(C_1), \dots, f(C_\beta)$$ with the following two properties:
\begin{description}
\item[(a) ] Every directed path $P_i~(i\in [\alpha])$ with positive flow connects a source vertex to a sink vertex.
\item[(b) ] $\alpha + \beta \le |V(D)| + |A(D)|$ and $\beta \le |A(D)|$.
\end{description}
\end{thm}

By Theorems~\ref{thm3-8} and~\ref{07-0}, Bai, Sun, Wang and Yu obtained the following complexity result for Eulerian digraphs (where $k$ is fixed and $\ell$ is part of the input).
\begin{thm}\label{07-1}\cite{Bai-Sun-Wang-Yu}
  Let $k \geq 3$ be a fixed integer and $\ell$ be part of the input. Let $D$ be an Eulerian digraph and $S\subseteq V(D)$ with $|S|=k$. The problem of deciding whether $\lambda_S^c(D)\geq \ell$ is NP-complete.
\end{thm}






The problem of {\sc Planar Arc-disjoint Paths (with Two Demand Pairs on Outer Face)} is formulated as follows:
\vspace{2mm}

\noindent {\bf Planar Arc-disjoint Paths (with Two Demand Pairs on Outer Face)}:
Let $D=(V,A)$ be a planar digraph and $(s_1,t_1), (s_2, t_2)$ be two demand pairs such that all terminals of the demand pairs lie on the outer face of a fixed planar embedding of $D$, appearing in the cyclic order $s_1,s_2,t_1,t_2$ along the outer face. The goal is to find $d_1$ $s_1-t_1$ paths, and $d_2$ $s_2-t_2$ paths such that all of the $d_1+d_2$ paths are arc-disjoint.

\vspace{2mm}

In 2012, Naves~\cite{Naves} established the NP-completeness for this problem.

\begin{thm}\label{07-2}~\cite{Naves}
The problem of {\sc Planar Arc-disjoint Paths (with Two Demand Pairs on Outer Face)} is NP-complete.
\end{thm}

By Theorem~\ref{07-2}, Bai, Sun, Wang and Yu proved the complexity for $\lambda_{S}^{c}(D)$ on planar digraphs (where $k$ is fixed and $\ell$ is part of the input).

\begin{thm}\label{07-3}\cite{Bai-Sun-Wang-Yu}
Let $k\geq 2$ be a fixed integer, $\ell $ be a part of the input. Let $D$ be a planar digraph and $S\subseteq V(D)$ with $|S|=k$. The problem of deciding whether \(\lambda _S^c(D) \geq \ell\) is NP-complete.
\end{thm}

The {\sc Hamiltonian Cycle} problem is a classic problem in graph theory. It is defined as follows:

\vspace{2mm}

\noindent {\bf Hamiltonian Cycle}: Given a graph $G$, the aim is to find a Hamiltonian cycle in $G$.

\vspace{2mm}
It is known that the Hamiltonian cycle problem is NP-complete for general undirected graphs, planar graphs and  $d$-regular graphs where $d \geq3$ ~\cite{M. R. Garey, Karp, C. Picouleau}.


\begin{thm}\label{07-4}~\cite{Karp}
The {\sc Hamiltonian Cycle} problem in undirected graphs is NP-complete.
\end{thm}

\begin{thm}\label{07-PlaHamil}\cite{M. R. Garey}
The {\sc Hamiltonian Cycle} problem in planar graphs is NP-complete.
\end{thm}

\begin{thm}\label{07-regular}~\cite{C. Picouleau}
For any fixed $d \geq 3$, the {\sc Hamiltonian Cycle} problem in $d$-regular graphs is NP-complete. 
\end{thm}

By Theorem~\ref{07-4}, Bai, Sun, Wang and Yu determined the complexity of $\lambda_S^c(D)$ on symmetric digraphs (where $k$ is part of the input and $\ell$ is fixed). 

\begin{thm}\label{07-5}\cite{Bai-Sun-Wang-Yu}
Let $\ell \geq 1$ be a fixed integer and $k$ be part of the input. Let $D$ be a symmetric digraph and $S\subseteq V(D)$ with $|S|=k$. The problem of deciding whether \(\lambda _S^c(D) \geq \ell\) is NP-complete.
\end{thm}

However, when $k\geq 2$ is fixed and $\ell=2$, the above problem in symmetric digraphs becomes polynomial-time solvable. This result and Theorem~\ref{07-5} imply the entries of Table~7-4.

\begin{thm}\label{07-6}\cite{Bai-Sun-Wang-Yu}
Let $k \geq 2$ be a fixed integer. Let $D$ be a symmetric digraph and $S\subseteq V(D)$ with $|S|=k$. The problem of deciding whether \(\lambda _S^c(D) \geq 2\) can be solved in polynomial time.
\end{thm}


By Theorem~\ref{07-PlaHamil}, Bai, Sun, Wang and Yu got the complexity for $\lambda _S^c(D)$ on planar digraphs (where $\ell$ is fixed and $k$ is part of the input). This result and Theorem~\ref{07-3} imply the entries of Table~7-5.

\begin{thm}\label{07-7}\cite{Bai-Sun-Wang-Yu}
Let $\ell \geq 1$ be a fixed integer and $k$ be part of the input. Let $D$ be a planar digraph and $S\subseteq V(D)$ with $|S|=k$. The problem of deciding whether \(\lambda _S^c(D) \geq \ell\) is NP-complete.
\end{thm}

By Theorem~\ref{07-regular}, Bai, Sun, Wang and Yu got the complexity for $\lambda _S^c(D)$ on Eulerian digraphs (where $\ell$ is fixed and $k$ is part of the input). This result and Theorem~\ref{07-1} imply the entries of Table 7-6.

\begin{thm}\label{07-8}\cite{Bai-Sun-Wang-Yu}
Let $\ell \geq 1$ be a fixed integer and $k$ be a part of the input. Let $D$ be an Eulerian digraph and $S\subseteq V(D)$ with $|S|=k$. The problem of deciding whether \(\lambda _S^c(D) \geq \ell\) is NP-complete.
\end{thm}

\begin{center}
	\begin{tabular}{|l|c|c|} \hline
		\multicolumn{3}{|c|}{Table 7-4: Symmetric digraphs} \\ \hline
		\makecell[l]{$\lambda_S^{c}(D) \geq \ell?$ $|S|=k$} & $k \geq 2$ constant & $k$ part of input\\
		\hline \hline
		$\ell =1,2$ & Polynomial\cite{Bai-Sun-Wang-Yu} & NP-complete\cite{Bai-Sun-Wang-Yu} \\ \hline
		$\ell \geq 3$ constant & & NP-complete\cite{Bai-Sun-Wang-Yu} \\ \hline
		$\ell$ part of input & & NP-complete\cite{Bai-Sun-Wang-Yu} \\ \hline
	\end{tabular}
\end{center}

\begin{center}
	\begin{tabular}{|l|c|c|} \hline
		\multicolumn{3}{|c|}{Table 7-5: Planar digraphs} \\ \hline
		\makecell[l]{$\lambda_S^{c}(D) \geq \ell?$ $|S|=k$} & $k \geq 2$ constant & $k$ part of input\\
		\hline \hline
		$\ell \geq 1$ constant & & NP-complete\cite{Bai-Sun-Wang-Yu} \\ \hline
		$\ell$ part of input & NP-complete\cite{Bai-Sun-Wang-Yu} & NP-complete\cite{Bai-Sun-Wang-Yu} \\ \hline
	\end{tabular}
\end{center}

\begin{center}
	\begin{tabular}{|l|c|c|c|} \hline
		\multicolumn{4}{|c|}{Table 7-6: Eulerian digraphs} \\ \hline
		\makecell[l]{$\lambda_S^c(D) \geq \ell?$ \\ $|S|=k$} & $k=2$ & $k \geq 3$ constant & $k$ part of input\\
		\hline \hline
		$\ell \geq 1$ constant & & & NP-complete\cite{Bai-Sun-Wang-Yu} \\ \hline
		$\ell$ part of input & & NP-complete\cite{Bai-Sun-Wang-Yu} & NP-complete\cite{Bai-Sun-Wang-Yu} \\ \hline
	\end{tabular}
\end{center}

\begin{op}\label{op5-0-3}
Complete the blanks in Tables~7-4, 7-5 and 7-6.
\end{op}

It would also be interesting to consider other classes of digraphs, such as semicomplete digraphs.
\begin{op}\label{op5-2}
Let $k\geq 2$ and $\ell \geq 1$ be fixed integers. What is the complexity of deciding whether $\lambda^c_{S}(D) \geq \ell$ for a semicomplete digraph $D$? where $S \subseteq V(D)$ and $|S|=k$.
\end{op}

K\"{u}hn and Osthus gave the following sufficient condition for a digraph to be $k$-cyclic.

\begin{thm}\label{thm5-5}\cite{Kuhn2008}
Let $k\geq 2$ be an integer. Every digraph $D$ of order $n\geq 200k^3$ which satisfies $\delta^0(D)\geq (n+k)/2-1$ is $k$-cyclic.
\end{thm}

Using Theorem~\ref{thm5-5}, we can prove the following:

\begin{pro}\label{pro-cycle}\cite{Sun4}
Let $k\geq 2, \ell\geq 1$ be integers. Let $D$ be a digraph of order $n\geq 200k^3$ which satisfies $\delta^0(D)\geq (n+k)/2+\ell-2$. Then $\lambda^c_S(D)\geq \ell$ for every $S\subseteq V(D)$ with $|S|=k$.
\end{pro}
\begin{pf}
Let $D$ and $S$ be defined as in the assumption. At stage 1, we set $D_0=D$. Since $\delta^0(D_0)\geq (n+k)/2+\ell-2\geq (n+k)/2-1$, by Theorem~\ref{thm5-5}, there is an $S$-cycle, say $C_1$, in $D_0$, then we obtain $D_1$ from $D_0$ by deleting the arcs of $C_1$. Generally, at stage $i\in [\ell]$, we start with a digraph $D_{i-1}$ which is obtained from $D_0$ by deleting the arcs of $i-1$ arc-disjoint cycles: $C_j$~$(j\in [i-1])$. Since now $\delta^0(D_{i-1})\geq (n+k)/2+\ell-2-(i-1)=(n+k)/2+\ell-i-1\geq (n+k)/2-1$, by Theorem~\ref{thm5-5}, there is an $S$-cycle, say $C_{i}$, in $D_0$, then we obtain $D_i$ from $D_{i-1}$ by deleting the arcs of $C_{i}$. By induction, we can obtain a set of $\ell$ arc-disjoint $S$-cycles: $C_j$~$(j\in [\ell])$. Therefore, $\lambda^c_S(D)\geq \ell$.
\end{pf}

\subsection{A related topic: directed cycle connectivity} 

By replacing the ``$S$-strong subgraphs" in the definition of strong subgraph $k$-connectivity \index{strong subgraph $k$-connectivity} with ``$S$-cycles", Wang and Sun~\cite{Wang-Sun} defined the following concept of {\em directed cycle $k$-connectivity} $\kappa _{k}^{c}(D)$ of a digraph $D$. Recall that $\kappa^c_S(D)$ denotes the maximum number of pairwise internally-disjoint $S$-cycles in $D$. Let 
$$\kappa _{k}^{c}(D)=\min\left \{ \kappa _{S}^{c}(D)\mid S\subseteq V(D),\left |S \right | =k,2\le k\le n \right \}.$$

In~\cite{Wang-Sun}, Wang and Sun studied the directed cycle $k$-connectivity of complete digraphs $\overleftrightarrow{K}_{n}$ and  complete regular bipartite digraphs $\overleftrightarrow{K}_{a,a}$. They gave a sharp lower bound for $\kappa _{k}^{c}(\overleftrightarrow{K}_{n})$ and determine the exact values for $\kappa _{k}^{c}(\overleftrightarrow{K}_{n})$ when $k\in \{2,3,4,6\}$. They also got the exact value of $\kappa _{k}^{c}(\overleftrightarrow{K}_{a,a})$ for each $2\leq k\leq n$.

After that, Bai, Sun, Wang and Yu~\cite{Bai-Sun-Wang-Yu} defined the following {\em directed cycle $k$-arc-connectivity} $\lambda _{k}^{c}(D)$ of a digraph $D$. Recall that $\lambda^c_S(D)$ denotes the maximum number of pairwise arc-disjoint $S$-cycles in $D$. Let 
$$\lambda _{k}^{c}(D)=\min\left \{ \lambda _{S}^{c}(D)\mid S\subseteq V(D),\left |S \right | =k,2\le k\le n \right \}.$$
	
In~\cite{Bai-Sun-Wang-Yu}, Bai, Sun, Wang and Yu determined the precise values of the directed cycle $k$-arc-connectivity of several digraph classes, including complete digraphs, complete bipartite digraphs and complete regular multipartite digraphs.

\vskip 1cm

\noindent {\bf Acknowledgement.} This work was supported by National Natural Science Foundation of China under Grant No. 12371352 and Yongjiang Talent Introduction Programme of Ningbo under Grant No. 2021B-011-G.


\end{document}